\documentclass[12pt,reqno]{article}
\usepackage{amsthm}
\usepackage{amsfonts}
\usepackage{amssymb}
\usepackage{amsmath}
\usepackage{epsfig}
\usepackage{graphics}

\newcommand{\R}{ { \mathbb{R} } }
\newcommand{\C}{ { \mathbb{C} } }
\newcommand{\T}{ { \mathbb{T}  } }
\newcommand{\N}{ { \mathbb{N} } }
\newcommand{\Q}{ { \mathbb{Q} } }
\newcommand{\Z}{ { \mathbb{Z} } }
\newcommand{\nz}{ { \setminus \{ 0 \} } }

\newcommand{\ep}{ \varepsilon }
\newcommand{\rep}{ \sqrt{\varepsilon} }

\def\ran{\operatorname{Ran} }
\newcommand{\Tr}{ \mbox{\rm Tr \!} }
\def\curl{\operatorname{curl}}
\def\Div{\operatorname{div}}

\newcommand{\dx}{ \partial_x }
\newcommand{\dy}{ \partial_y }
\newcommand{\Dy}{ {\rm D}_y }
\newcommand{\dz}{ \partial_z }
\newcommand{\Dz}{ {\rm D}_z }
\newcommand{\dt}{ \partial_t }

\newcommand{\dT}{ \partial_T }

\newcommand{\Dth}{ {\rm D}_\theta }
\newcommand{\dthu}{ \partial_{\theta_1} }
\newcommand{\dthz}{ \partial_{\theta_0} }

\newcommand{\tendlorsque}[2]{\mathop{\longrightarrow}\limits_{#1\rightarrow#2}}

\newcommand{\fin}{ \hfill $\square$ \vspace{0.2cm}  }
\newcommand{\ie}{\emph{i.e.} }
\newcommand{\cf}{\emph{cf.} }
\newcommand{\via}{\emph{via} }

\def\Og{\operatorname{\Omega_\gamma}}
\def\aO{\operatorname{ad_\Omega}}
\newcommand{\Lep}{ L\left(\dt,\dx,\frac{1}{\rep}\dy,\frac{1}{\rep}\dz\right) }
\newcommand{\Lsing}{ L\left(\dt,\dx+\frac{1}{\ep}\dthz,\frac{1}{\rep}\dy,
\frac{1}{\rep}\dz\right) }
\newcommand{\crd}[2]{[#1,#2]_{\rm d}}
\newcommand{\crod}[2]{[#1,#2]_{\rm od}}

\newcommand{\fiel}{ {\bf u} }
\newcommand{\pop}{ {\bf N} }
\newcommand{\coh}{ {\bf C} }
\newcommand{\sol}{ {\bf U} }
\newcommand{\fiele}{ {\bf u}^\ep }
\newcommand{\pope}{ {\bf N}^\ep }
\newcommand{\cohe}{ {\bf C}^\ep }
\newcommand{\sole}{ {\bf U}^\ep }
\newcommand{\solc}{ \mathcal{U}^\ep }
\newcommand{\fiela}{ {\bf u}^\ep_{\rm app} }

\newcommand{\coha}{ {\bf C}^\ep_{\rm app} }

\newcommand{\sola}{ {\bf U}^\ep_{\rm app} }
\newcommand{\solac}{ \mathcal{U}^\ep_{\rm app} }
\newcommand{\fzo}{ {\fiel_{\rm osc}^0} }

\newcommand{\fzz}{ {\fiel^{0,0}} }

\newcommand{\fuz}{ {\fiel^{1,0}} }
\newcommand{\fuzo}{ {\fiel_{\rm osc}^{1,0}} }
\newcommand{\fdo}{ {\fiel_{\rm osc}^2} }
\newcommand{\fdu}{ {\fiel^{2,1}} }
\newcommand{\fdzo}{ {\fiel^{2,0}_{\rm osc}} }
\newcommand{\fdz}{ {\fiel^{2,0}} }
\newcommand{\pzz}{ {\pop^{0,0}} }

\newcommand{\puz}{ {\pop^{1,0}} }
\newcommand{\puu}{ {\pop^{1,1}} }

\newcommand{\czu}{ {\coh^{0,1}} }
\newcommand{\cuz}{ {\coh^{1,0}} }
\newcommand{\cuu}{ {\coh^{1,1}} }
\newcommand{\dep}{ {\delta^\ep} }
\newcommand{\Dep}{ {\Delta^\ep} }

\newcommand{\vo}{k}
\newcommand{\vf}{k}

\newcommand{\ch}{ \mathcal{C} }
\newcommand{\cp}{ \mathcal{C}_+ }
\newcommand{\cm}{ \mathcal{C}_- }
\newcommand{\cz}{ \mathcal{C}_0 }
\newcommand{\res}{ \mathcal{R}(\vo) }

\def\roteps{\operatorname{curl_\ep}}
\newcommand{\op}{ { \, \sharp \, } }

\newcommand{\var}{ (t,x,y,z) }
\newcommand{\varT}{ (t,x,y,z,T) }
\newcommand{\varTt}{ (t,x,y,z,T,\theta) }
\newcommand{\varTdt}{ (t,x,y,z,T,\sigma,\theta) }

\newcommand{\ea}{ ^\ep_{\rm app} }

\newcommand{\mun}{ M_1(\dt,\dx) \,}

\newcommand{\muda}{ 
M_1(-\gamma\kappa-i\vf\cdot\alpha_1,i\vo\cdot\alpha_0) \,}
\newcommand{\muth}{
  M_1(-\vf\cdot\partial_{\theta_1},\vo\cdot\partial_{\theta_0}) \,} 
\newcommand{\mudth}{ M_1(-\gamma\kappa-
\vf\cdot\partial_{\theta_1},\vo\cdot\partial_{\theta_0}) \,} 
\newcommand{\muai}{ M_1(-\vf\cdot\alpha_1,\vo\cdot\alpha_0)^{-1} \,}

\newcommand{\muthi}{
  M_1(-\vf\cdot\partial_{\theta_1},\vo\cdot\partial_{\theta_0})^{-1} \,} 
\newcommand{\mudthi}{
  M_1(-\gamma\kappa-\vf\cdot\partial_{\theta_1},\vo\cdot
\partial_{\theta_0})^{-1} \,} 
\newcommand{\md}{ M_2\,(\dT,\dy,\dz) \,} 
\newcommand{\mdyz}{ M_2\,(0,\dy,\dz) \,}  
\newcommand{\val}{ {\rm v}(\alpha) } 
\newcommand{\vDth}{ {\rm v}(\Dth) } 
\newcommand{\aal}{ {\rm a}(\alpha) } 
\newcommand{\aDth}{ {\rm a}(\Dth) } 
\newcommand{\pk}{ p_k(\Dy,\Dz) }

\newtheorem{lemma}{Lemma}[section]

\newtheorem{theo}[lemma]{Theorem}

\newtheorem{prop}[lemma]{Proposition}
\newtheorem{lemme}[lemma]{Lemma}
\newtheorem{hyp}[lemma]{Hypothesis}
\newtheorem{defn}[lemma]{Definition}
\newtheorem{rem}[lemma]{Remark}
\newtheorem{ex}[lemma]{Example}
\newtheorem{nota}[lemma]{Notation}

\title{High frequency behaviour of the Maxwell-Bloch model with relaxations: convergence to
 the Schr\"odinger-Boltzmann system} 
\author{F. Castella $^{(1)}$ and E. Dumas $^{(2)}$}
\date{}

\begin{document}

\maketitle

\begin{center}
\small
(1) IRMAR, UMR 6625 (CNRS-UR1) \\
Universit\'e de Rennes 1\\
Campus de Beaulieu, 35042 Rennes Cedex - France\\
email: francois.castella@univ-rennes1.fr
\end{center}

\begin{center}
\small
(2) Institut Fourier, UMR 5582 (CNRS-UJF) \\
100 rue des Math\'ematiques \\
Domaine Universitaire \\
BP 74, 38402 Saint Martin d'H\`eres - France \\
email: edumas@ujf-grenoble.fr
\end{center}

\begin{abstract}
We study the Maxwell-Bloch model, which describes the propagation of a laser through a material and the associated interaction between laser and matter
(polarization of the atoms through light propagation, photon emission and absorption, etc.).
The laser field is described through Maxwell's equations, a classical equation,
while matter is 
represented at a quantum level and satisfies a quantum Liouville equation known as the Bloch model.
Coupling between laser and matter is described through
a quadratic source term in both equations.
The model also takes into account partial relaxation effects, namely the trend of matter to return to its natural thermodynamic equilibrium.
The whole system involves $6$+$N \, (N+1)/2$ unknowns, the six-dimensional electromagnetic field
plus the $N \, (N+1)/2$ unknowns
describing the state of matter, where $N$ is the number of atomic energy levels of the considered material.

We consider at once a high-frequency and weak coupling situation,
in the general case of anisotropic electromagnetic fields that are subject to 
diffraction. Degenerate energy levels are allowed.
The whole system is stiff and involves strong nonlinearities.

\medskip

We show the convergence to a nonstiff,
nonlinear, coupled
Schr\"o\-din\-ger-Boltzmann model, involving $3$+$N$ unknowns. The electromagnetic field is eventually
described through its envelope, one unknown vector in $\C^3$. It satisfies a Schr\"o\-din\-ger
equation that takes into account propagation and diffraction of light
inside the material. Matter on the other hand  is described through
a $N$-dimensional vector describing the occupation numbers of each atomic level. It satisfies a Boltzmann equation that describes
the jumps of the electrons between the various atomic energy levels,
as induced by the interaction with light.
The rate of exchange between the atomic levels is proportional to
the intensity of the laser field.
The whole system is the physically natural nonlinear model.

In order to provide an important and explicit example, we completely 
analyze the specific (two dimensional) 
Transverse Magnetic case, for which formulae turn out to be simpler.

Technically speaking, our analysis 
does not enter the usual mathematical framework of geometric 
optics: it is more singular, and requires 
an \emph{ad hoc} Ansatz.
\end{abstract}

\tableofcontents


\vspace{5mm}
\noindent
Support by the program ``SYDYQ: SYst\`emes DYnamiques Quantiques''  
from the Universit\'e Joseph Fourier (Grenoble) is acknowledged.

\section{Introduction}


Maxwell-Bloch systems are of common use in Laser Physics 
(see the textbooks \cite{Boh79}, \cite{Boy92}, 
\cite{CTDRG88}, \cite{NM92}, \cite{PP69}, \cite{SSL77}). 
They modelize the evolution of an electromagnetic field, treated classically,
and coupled with an ensemble of identical atoms, which in turn are described 
by a quantum density matrix. This model is relevant when
atoms are far from the ionization energy (to possess discrete energy levels), while
they have sufficiently low density and the laser field is strong enough
(which allows to describe the field classically while matter is described in a quantum way
 -- see \cite{CTDRG88}). 

In the Maxwell-Bloch model, the electromagnetic field satisfies
Maxwell's equations, whose unknowns are
the electric and magnetic fields $E\in \R^3$ and $B\in \R^3$,
where $E=E(t,x,y,z)$ and $B=B(t,x,y,z)$ and 
$t \in \R$ is time while $(x,y,z)\in\R^3$ are the space coordinates.
Matter is described through a Bloch equation, whose unknown is
the density matrix $\rho=\rho(t,x,y,z)$, a quantum variable which describes
the atomic state at $(t,x,y,z)$. 
We consider that the atoms only visit
the $N$ lowest energy levels. 
The latter are the $N$ first eigenstates of
the free material system, in the absence of 
fields. In this basis, the density matrix $\rho=\rho(t,x,y,z)$
is an $N\times N$ matrix, for each value of $(t,x,y,z)$.
The diagonal entries 
$\rho(t,x,y)(n,n)$ (called \emph{populations}) give the 
proportion of matter that lies in the $n$-th energy level ($n=1,\ldots,N$), while the
off-diagonal entries $\rho(t,x,y,z)(n,p)$ with $n \neq p$, (called \emph{coherences}) 
give the correlation between levels $n$ and $p$. 
The complete Maxwell-Bloch system takes into account
the coupling between the laser field and the atoms \via terms  that are quadratic, proportional to $\rho \times E$,
and which describe polarization of matter due to laser propagation.

\medskip

We study the high frequency and weak coupling behaviour
of the Maxwell-Bloch system, a situation in which the typical frequencies of the field and of
the atoms' oscillations are large and possibly resonate, while the strength of laser-matter coupling is small.

For large frequencies, the electromagnetic field 
is expected to be asymptotically solution to a nonlinear Schr\"odinger 
equation. This is the paraxial approximation. We refer to \cite{DJMR95} and 
\cite{JMR98}, as well as \cite{Lan98} on these matters, when the sole laser field
propagates (no coupling with matter). When matter is actually coupled to the field,
we refer to \cite{BBCNT04}. Here
a high frequency Maxwell-Bloch system is studied both physically and mathematically, for atoms that only possess $3$ non-degenerate energy 
levels. The analysis leads to a Schr\"odinger-Bloch approximation of the original
system, in a spirit similar to the present paper.

The weak coupling behaviour of matter is a bit more delicate
to handle: to have a clean limit, one needs
to take thermodynamic fluctuations into account. For this reason we introduce,
in a standard fashion,
phenomenological \emph{relaxation operators} in the original Bloch system.
These impose 
a rapid decay of coherences, as well as a quick return to  equilibrium
of populations.
We refer to \cite{BBR01} for 
mathematical properties of the relaxation operators that are natural in this context. 
Due to the relaxation effects, it is expected that 
Bloch's equation is asymptotic to a Boltzmann equation, 
sometimes called ``Einstein's rate equation'' (\cf 
\cite{Lou91}, \cite{BB06}). It describes how the atoms jump
between the various energy levels under the action of the external field.
When the driving high frequency field is given (and the Bloch model is thus linear),
we refer to the papers \cite{BFCD04} and \cite{BFCDG04}, which  study the actual convergence 
of Bloch's equation to a Boltzmann model in the 
weak coupling regime.
In that case, a formula is found for the transition rates
involved in the limiting Einstein rate equation, which coincides
with the one formally obtained in the Physics literature.
Note however that the question studied in \cite{BFCD04} and \cite{BFCDG04}
is a linear problem, and proofs strongly  use ODE averaging techniques
as well as the positivity of relaxation operators (features that the present
text does not share).
We also mention \cite{CDG07}, where similar asymptotics are 
treated both for quantum 
and classical models.
We stress finally that many other works deal with the rigorous derivation of 
Boltzmann like equations from (usually linear) models 
describing the interaction of waves/particles with external media. 
A non-convergence result is given in \cite{CP02} and 
\cite{CP03}. Convergence in the case of an electron 
in a periodic box is studied in \cite{Cas99}, 
\cite{Cas01}, \cite{Cas02}, while the case of an 
electron in a random medium is addressed in 
\cite{EY00}, \cite{KPR96}, \cite{Spo77}, \cite{Spo80}, 
\cite{Spo91} -- see also \cite{Nie96} for a semi-classical approach. In a nonlinear context, a partial 
result is obtained in \cite{BCEP04}.

\medskip

The above  formal discussion suggests, in the present case, that
the high-frequency Maxwell system goes to a Schr\"odinger model for the
envelope of the field, while the weakly coupled Bloch system supposedly goes to a Boltzmann equation describing the jumps of electrons
between the atomic levels. 

This is the program we rigorously develop in the present paper. We fully prove convergence 
of the coupled Maxwell-Bloch system to a coupled Schr\"odinger-Boltzmann model. We also prove that 
the rate of exchange between the energy levels
is proportional to the laser's intensity.
In doing so we recover the physically relevant model.
Our approach mainly uses three-scales geometric optics, yet in a more singular
context where the partial
relaxation effects impose a specific treatment of coherences.


\section{Presentation of the results}


\subsection{The model}

The Maxwell-Bloch system, whose unknowns are the electric field
$E=E(t,x,y,z)\in \R^3$, the magnetic field $B=B(t,x,y,z)\in\R^3$, together
with the density matrix $\rho=\rho(t,x,y,z) \in \C^{N\times N}$ (the
space of $N\times N$ complex matrices),
reads
\begin{align}
\label{Max1} 
&
 \dt B + \curl E = 0 , \\
&
\label{Max2}
 \dt E - \curl B = -\dt P , \qquad \text{ with } \, P=\Tr(\Gamma\rho),
\\
&
\label{Bloch}
\dt \rho = -i [\Omega - E \cdot \Gamma , \rho] + Q(\rho),
\\\
&
\text{where } \,
\label{qrho}
Q(\rho)=W\op\rho_{\rm d} - \gamma \, \rho_{\rm od}.
\end{align}
In these equations, $\curl$ is the usual $\curl$ operator on vector fields in $\R^3$,
$\Gamma$ as well as $\Omega$ are given matrices in $\C^{N\times N}$, $\gamma>0$ is
a given positive
constant, and whenever $A$ and $B$
are matrices in $\C^{N\times N}$, the trace  $\Tr(A)$ denotes the usual trace of $A$ while the
bracket $[A,B]$ denotes the usual commutator between matrices
$$[A,B]=AB-BA.$$
The term $Q(\rho)=Q(\rho)(t,x,y,z)$ is the so-called relaxation matrix ,
an $N \times N$ matrix for each value of $(t,x,y,z)$.
Its definition involves 
$W \in \C^{N\times N}$, a given 
matrix with nonnegative entries $W(n,m)\geq 0$, while $\rho_{\rm d}$ and 
$\rho_{\rm od}$ denote
the diagonal respectively off-diagonal parts of the 
density matrix $\rho$ (they correspond respectively to the populations and the coherences).
They are $N \times N$ matrices defined,
for each value of $(t,x,y,z)$, by their entries
\begin{align*}
&
\rho_{\rm d}(t,x,y,z)(n,p)=\rho(t,x,y,z)(n,p) \, {\bf 1}[n=p],
\\
&
\rho_{\rm od}(t,x,y,z)(n,p)=\rho(t,x,y,z)(n,p) \, {\bf 1}[n\neq p].
\end{align*}
Equation \eqref{qrho} also uses the following notation,
valid thoughout the present text : given any matrix $A$ with nonnegative entries we set
\begin{align}
\label{wdiese}
\left\{
\begin{array}{l}
\vspace{0.1cm}
A\op \rho_{\rm d}(n,n) = 
\displaystyle
\sum_{k=1}^N \left[ A(k,n) \, \rho_{\rm d}(k,k)-A(n,k) \, \rho_{\rm d}(n,n)\right],
\\
A\op \rho_{\rm d}(n,p) =0 \quad \text{ when } \, n\neq p.
\end{array}
\right.
\end{align}
The meaning of operator $Q(\rho)$ in \eqref{Bloch} is the following. The term
$-\gamma \, \rho_{\rm od}$
induces exponential relaxation  to zero for the coherences, while
the term $W\op\rho_{\rm d}$ acts on the populations only, and induces exponential relaxation
of the populations towards some thermodynamical equilibrium that depends on the values of the
$W(n,p)$'s. As in conventional kinetic theory of gases,
relation \eqref{wdiese} asserts that along time evolution, atoms may leave with probability $W(k,n)$ 
the $k$th eigenstate
to populate the $n$th eigenstate (this is the so-called gain term
$\sum_{k=1}^N W(k,n) \, \rho_{\rm d}(k,k)$), while some atoms may conversely leave with probability $W(n,k)$
the $n$th state
to populate some other $k$th state
(this is the so-called loss term
$- \sum_{k=1}^N W(n,k) \, \rho_{\rm d}(n,n)$).

Theoretically, the Maxwell-Bloch system
needs to be supplemented with the Amp\`ere and Faraday laws, 
\begin{equation} \label{far} 
\Div B = 0, \qquad \Div (E+P) = 0.
\end{equation} 
These constraints (\ref{far}) are anyhow transported as soon as they are satisfied by the
initial data, hence we shall skip them in the sequel.

Note that the above equations are readily given in the convenient
dimensionless form that suits our purpose. The precise scaling under study is discussed later.

\medskip

We now comment on these equations, and on all involved quantities.

The density matrix $\rho(t,x,y,z)$ is Hermitian and  
positive. It describes the state of matter at $(t,x,y,z)$.

The constant matrix $\Omega$ is the
free Hamiltonian of the material system, written in the natural
eigenbasis. It is a fixed physical constant associated with the considered atomic species. It reads
\begin{align}
\label{OOMM}
\Omega = {\rm diag}\,(\omega(1),\dots,\omega(N)),
\end{align}
where $0<\omega(1)\leq\cdots\leq\omega(N)$ are the atomic energies. For later convenience, we readily introduce
the differences between energy levels, as
\begin{align}
\label{omnk}
\omega(n,k)=\omega(n)-\omega(k).
\end{align}

The constant matrix $\Gamma$ is called the dipolar operator.
It is a hermitian matrix.
It has the value 
$N_{\rm atomic} \times \gamma_{\rm atomic}$,
where $N_{\rm atomic}$ is the number of atoms per unit volume, while
$\gamma_{\rm atomic}$ is a fixed physical constant (a matrix) associated with the considered atomic species.

The entries of so-called dipolar momentum $E\cdot\Gamma$ 
are defined for any $m,n =1, \ldots ,N$,  as
$(E\cdot\Gamma)(t,x,y,z)(m,n)= E(t,x,y,z)\cdot\Gamma(m,n)$,
where $a \cdot b$ denotes the componentwise product of two vectors in $\C^3$, namely
$$
(E\cdot\Gamma)(t,x,y,z)(m,n) 
= E_x(t,x,y,z) \Gamma(m,n)_x + E_y(t,x,y,z) \Gamma(m,n)_y
+ E_z(t,x,y,z) \Gamma(m,n)_z, 
$$
and the subscripts $x$ (resp. $y$, resp. $z$), denote the $x$ (resp. $y$, resp. $z$) components of the relevant vectors
(note the absence of complex conjugation).

The non-negative transition coefficients $W(n,k)\geq 0$ are  known as
the Pauli 
coefficients, see \cite{Boy92}.
They
satisfy a micro-reversibility relation at temperature $T$, {\em i.e.}
$W(n,k) = W(k,n) \exp\left(\omega(k,n)/T\right)$. In that perspective the relaxation
operator $Q(\rho)$ translates the fact that the atoms tend to relax towards the thermodynamical equilibrium given by 
$\rho_{\rm od}\equiv 0$ and $\rho_{\rm d}(n,n)\equiv\exp\left(-\omega(n)/T\right)$.
Both the off-diagonal relaxation term $\gamma$ and the Pauli coefficients $W(n,k)$ are
physical data. Contrary to $\Gamma$ or $\Omega$, their relation with given physical
constants attached with the specific atomic species at hand is unclear.

The Bloch equation (\ref{Bloch}) relies on the so-called dipolar approximation : 
the coupling between light and matter is taken into account through the simplest
$[E\cdot\Gamma,\rho]$ term, which is quadratic (proportional to $\rho\times E$) and local (it only depends
on the value of $E$ and $\rho$ and the {\em same} point $(t,x,y,z)$). This approximation
implicitely assumes that the wavelength of the field is larger than the typical spatial extension of the atom, hence can be taken constant over the whole domain occupied by each given atom.

\subsection{The scaling}

Let us introduce the physical scales in the model, which transform all
constants ($\Gamma$, $\Omega$, $\gamma$, $W$) and unknowns ($E$, $B$, $\rho$)
into quantities of order one in the regime we wish to study.

\medskip

Firstly, concerning the time variations, we study a high frequency regime.
Calling $1/\ep$ the dimensionless parameter measuring typical
values of the frequencies, {\em i.e.} the ratio between the time scale of observations
and the time scale of the variations of $E$, $B$ and $\rho$,
all time derivatives $\dt$ then become $\ep\dt$ in the scaled version of \eqref{Max1}-\eqref{qrho}.

Secondly, we want to study a situation where
constructive interference occurs between the time oscillations of the electromagnetic field
and the ones of the atom  (giving rise to atomic absorption and emission of photons).
For that reason, the time
variations of $E$, $B$, and $\rho$ should all take place at similar frequencies. For that reason
$\Omega$ is naturally a  quantity of order one in the scaled model.

Thirdly, the typical strength of the coupling between light and matter is entirely 
determined by the physical constant $\Gamma$. 
Since weak coupling is realized when the polarization 
operator has an effect of the order one on the chosen time scale,
it turns out that we need to prescribe $\Gamma = \mathcal{O}(\rep)$.
Mathematically, this means we shall replace $\Gamma$ by $\rep \, \Gamma$ in the scaled model.
Indeed, the electric dipole momentum then becomes  
$|E\cdot(\rep \, \Gamma)|^2=\mathcal{O}(\ep)$ which, integrated over macroscopic time scales 
$\mathcal{O}(1/\ep)$, results in an energy 
contribution of order one as desired. Physically, since
$\Gamma=N_{\rm atomic} \times \gamma_{\rm atomic}$ where $N_{\rm atomic}$ is the
number of atoms per unit volume and $\gamma_{\rm atomic}$ is the polarizability
of one atom, this means
that the atomic density is here {\em tuned} to be of order $\mathcal{O}(\sqrt{\ep})$, so as to observe order one effect of coupling over this scale.

Fourthly, concerning the space variables, the hyperbolic nature of
the Max\-well equations \eqref{Max1}-\eqref{Max2}
suggests to rescale space so that $\dx$, $\dy$, $\dz$
become $\ep \dx$, $\ep \dy$, and $\ep \dz$, respectively, due to the scaling $\dt \mapsto \ep \dt$
and to finite propagation speed: the laser visits space scales $\mathcal{O}(1/\ep)$
over time scales  $\mathcal{O}(1/\ep)$.
This simple scaling would in fact be easily described along lines similar to the present analysis.
We wish to investigate physically richer situations where diffraction occurs. 
To this end, as in \cite{BBCNT04}, \cite{Dum03} and \cite{Dum04}, we impose 
anisotropy and introduce a third scale: on the one hand, we choose one direction of propagation, say the $x$ direction, meaning that at the microscopic scale, the fields vary with $x$ but not with $(y,z)$;
on the other hand, we restrict our attention to fields that slowly vary
over the scale $1/\ep$ in the $x$ direction (or, in other words, that have spatial extension
$1/\ep$ in $x$), while they slowly vary over the scale $1/\rep$ in the $(y,z)$ direction (or, in other words, that have spatial extension $1/\rep$ in $(y,z)$).
The typical shape of the laser beam is  thus that of a ``light cigar'', as in \cite{Don94}. 
In macroscopic scales, this provides fields that strongly oscillate at frequency $1/\ep$
in $x$ (and only in this direction),  while they have
support of size $1$ in $x$, and support of size $\rep$ in $(y,z)$. 
All this imposes to rescale $\dx$ as $\ep\, \dx$, and $(\dy,\dz)$ as
$(\rep \, \dy, \, \rep \, \dz)$ 
in the above equations.
We stress that this is a definite choice of shape of the kind of laser beams we wish to study.
It is also a choice of polarization: these beams are shot in the $x$ direction only.

Fifthly, consider relaxations $Q(\rho)$. In order for the 
diagonal relaxation $W\op\rho_{\rm d}$ to have an $\mathcal{O}(1)$ effect 
at times $t$ of order $\mathcal{O}(1)$, we take the coefficients in $W$ of size 
$\ep$, and write $\ep W(n,k)$ instead of $W(n,k)$. The 
off-diagonal relaxation $-\gamma \rho_{\rm od}$ on the other hand is supposed to have a much shorter time 
scale, which we choose to be of order $\mathcal{O}(\ep)$. This corresponds to 
an off-diagonal relaxation $-\gamma \, \rho_{\rm od}$ that remains unscaled. In fact, there 
is no theoretical description of the relaxation time as a 
function of $\ep$. In \cite{BFCD04} and \cite{BFCDG04}, the 
off-diagonal relaxation is scaled as $-\gamma\ep^\mu\rho_{\rm od}$ 
with an extra free parameter $\mu$ that is constrained to satisfy $0\leq\mu<1$ (to have a clean limit,
off-diagonal relaxation should be strong enough with respect to
the chosen time scale). For technical 
reasons, in these papers, $\mu$ is actually restricted to 
$0\leq\mu<1/4$. In the present paper, the formal analysis could be 
performed for $0\leq\mu<1$, yet our main stability result 
(Theorem~\ref{ThStab}) requires the stronger constraint $\mu=0$.

\medskip

After rescaling all variables and physical constants accordingly,
the system \eqref{Max1}-\eqref{qrho} becomes 
\begin{align}
\label{MBeps} 
\left\{
\begin{array}{l}
\vspace{0.2cm}
\dt B^\ep + \roteps E^\ep = 0,
\\
\vspace{0.2cm}
\dt E^\ep - \roteps B^\ep = 
\displaystyle\frac{i}{\rep} \, \Tr \left(\Gamma\Og\cohe\right) 
- i \Tr \left(\Gamma[E^\ep\cdot\Gamma,\cohe+\pope]\right) 
- \rep \, \Tr \left(\Gamma \, \,  W\op \pope\right),
\\
\vspace{0.2cm}
\dt\cohe = -\displaystyle\frac{i}{\ep} \Og \cohe
+ \frac{i}{\rep} \, \crod{E^\ep\cdot\Gamma}{\cohe+\pope},
\\
\vspace{0.2cm}
\dt\pope = \displaystyle\frac{i}{\rep} \, \crd{E^\ep\cdot\Gamma}{\cohe} 
+ W\op\pope.
\end{array}
\right.
\end{align}
Here and in the sequel, we adopt for convenience the  notation $\cohe$ and 
$\pope$ for coherences $\rho_{\rm od}$ and populations $\rho_{\rm d}$, respectively.
Besides, throughout the sequel, 
subscripts ``d'' and ``od'' shall always refer to diagonal and 
off-diagonal parts of the considered $N \times N$ matrices. Lastly, we also denote by $\roteps$ 
the curl operator associated with our scaling, {\em i.e.}
\begin{align}
\label{curle}
\roteps E = 
\left( \frac{1}{\rep}\dy E_z-\frac{1}{\rep}\dz E_y,
\frac{1}{\rep}\dz E_x-\dx E_z,\dx E_y-\frac{1}{\rep}\dy E_x
\right),
\end{align}
and $\Og$ is a shorthand notation for $\aO-i\gamma$, 
\begin{equation} \label{Og}
\Og\coh=[\Omega,\coh]-i\gamma\coh, \quad
\text{\em i.e.  }
\left(\Og\coh\right)(n,p)=\omega(n,p)\coh(n,p)-i\gamma\coh(n,p).
\end{equation}

\subsection{Description of the results}

\paragraph{Main result: profiles, separation of scales, and obtention of an approximate solution.}

Maxwell-Bloch's system~\eqref{MBeps} is a nonlinear hyperbolic 
symmetric system, singular in $\ep$, that we may write 
symbolically
$$\Lep\sole=F^\ep(\sole).$$
Its unknown is
$$\sole=(\fiele,\cohe,\pope),$$
{where } $\fiele$ \text{ stands for the sole electromagnetic field }
$$\fiele=(B^\ep,E^\ep).$$

We are interested in solving a Cauchy 
problem associated with \eqref{MBeps} for initial data that are smooth, but high 
frequency,  of order $\mathcal{O}(1/\ep)$. 
The typical 
difficulty is to ensure existence of the whole family 
$(\sole)_{\ep\in]0,\ep_0]}$ on some time interval $[0,t_\star]$,
\emph{independent} of $\ep$.

This 
enters the framework of $3$-scales diffractive optics (see 
\cite{DJMR95}, \cite{Lan98}, and the surveys \cite{JMR99}, 
\cite{Dum06}).  To this end, our analysis uses 
a WKB analysis, based on profiles that we expand in successive powers of $\ep$. 
An important  original point yet is,
we consider amplitudes that are even higher (namely $\mathcal{O}(1)$) than the ones 
allowed by transparency properties (namely $\mathcal{O}(\ep)$ -- see \cite{JMR00}). 
We show these larger amplitudes are eventually compensated by the {\em partial} relaxation: 
remember the sole equation on the coherences $\cohe$ carries the relaxation term
$-\gamma  \cohe / \ep$ in (\ref{MBeps}).
We refer to remark \ref{RemAnsatz} below
on this important point.
\medskip

We start with a given a wavevector
$$
k=\left(\vo^1,\ldots,\vo^d\right)\in\R^d, 
$$
and initial data of the form
\begin{align} 
\label{wn}
\sole_{\rm ini}(x,y,z) =
\underline{\sol^0}\left(x,y,z,\frac{k x}{\ep}\right)
+ \delta^\ep\left(x,y,z,\frac{k x}{\ep}\right).
\end{align}
The profiles $\underline{\sol^0}(x,y,z,\theta_0)$, 
$\delta^\ep(x,y,z,\theta_0)\in\mathcal{C}^\infty(\R³\times\T^d)$ that are
periodic in $\theta_0\in\R^d$ and  we assume that 
$$\delta^\ep\tendlorsque{\ep}{0}0$$
in every Sobolev space 
$H^s$. We also assume that $\underline{\sol^0}$ is submitted to {\em polarization
conditions}, namely that the electromagnetic field $\fiel^0$ has the usual directional constraints, see (\ref{polarFiel}), 
and that  some preferred components 
$\underline{\coh^0}(m,n)$ vanish, see \eqref{polarCoh}.

We 
build for all $\ep>0$ an approximate solution $\sola$  to (\ref{MBeps}),  defined
on some time interval $[0,t_\star]$, and with initial values 
$\underline{\sol^0}(x,y,z,\vo x/\ep)$. Our construction uses a leading profile 
$\sol^0$ and  two correctors $\sol^1$, $\sol^2$, and $\sola$ is of the form
\begin{equation} \label{sola}
\sola(t,x,y,z) = \sum_{j=0}^2 \rep^j\sol^j
(t,x,y,z,T,\sigma,\theta_0, \theta_1)\Bigg|_{
T=t/\rep, \,
\sigma={\gamma t}/{\ep},(\theta_0, \theta_1)=\left({\vo x}/{\ep} , -{\vf t}/{\ep}\right)} .
\end{equation}

Correctors $\sol^1$, $\sol^2$, are introduced so as to ensure consistancy of 
this Ansatz, namely
\begin{equation} \label{consist}
\Lep\sola-F^\ep(\sola)=\mathcal{O}(\rep).
\end{equation}
One main point in our analysis is the separation of scales, and the 
crucial introduction of the fast variables $T$, $\sigma$ and $\theta_1$ in (\ref{sola}). Explanations 
on this point are postponed to the next paragraph.

Once the approximate solution is constructed, we prove 
{\bf Theorem~\ref{ThStab}.} It asserts that, for $\ep_0>0$ small enough, 
there is a unique solution $\sole$ to system \eqref{MBeps} 
with initial value $\sole_{\rm ini}$, which is well approximated by  $\sola$,
namely
\begin{equation} \label{stab}
\forall\mu\in\N^3, \quad 
\|\partial^\mu_{x,y,z}\left(\sole-\sola\right)\|
_{L^\infty([0,t_\star]\times\R^3)}\tendlorsque{\ep}{0}0 .
\end{equation}
On top of that, we are able to completely describe the dynamics of the dominant
term $\sol^0$. It provides an $o(1)$ approximation of the original dynamics of
$\sole$. The function $\sol^0$ satisfies a coupled, nonlinear, Schr\"odinger-Bloch system, which
we describe later.

This result  is achieved \via a singular system method 
(see \cite{JMR95}), where the unknown is a profile (thus with 
non-singular initial data). The difficulty comes from the 
``supersingular'' nature of the system, and the fact that 
relaxations are only \emph{partial} ones (acting on a part 
of the dependent variables only). This problem is overcome 
thanks to the structure of the approximate solution.

\paragraph{Describing the asymptotic dynamics (1) - oscillations and initial layer.}

In the chosen Ansatz \eqref{sola}, we assume that
each profile $\sol^j$ may be decomposed into modes 
$\sol_\alpha^{j,\kappa}$, as
\begin{align}
\label{momo1}
\sol^j\varTdt=\sum_{\alpha\in\Z^{2d}}\sum_{\kappa\in\N}
\sol_\alpha^{j,\kappa} \, \exp\left(i\,\alpha\cdot\theta\right) \, \exp\left(-\kappa\sigma\right).
\end{align}
This defines the quantities $\sol_\alpha^{j,\kappa}=\sol_\alpha^{j,\kappa}(t,x,y,z,T)$, for which $\alpha \in \Z^{2d}$
is seen as a Fourier mode, while $\kappa\in\N$ is seen as an exponentially decaying mode.
We introduce in passing and for later convenience the notation
\begin{align}
\label{momo2}
\sol_\alpha^{j}
=\sum_{\kappa\in\N}
\sol_\alpha^{j,\kappa} \, \exp\left(-\kappa\sigma\right),
\end{align}
a function of $(t,x,y,z,T,\sigma)$, as well as
\begin{align}
\label{momo3}
\sol^{j,\kappa}
=
\sum_{\alpha\in\Z^{2d}}
\sol_\alpha^{j,\kappa} \, \exp\left(i\,\alpha\cdot\theta\right),
\end{align}
a function of $(t,x,y,z,T,\theta)$.
For $\kappa>0$, the decomposition (\ref{momo1}) encodes an 
exponential mode in the
$\sigma=\gamma t/\ep$ variable,
representative of the initial layer induced by off-diagonal 
relaxations, while
Fourier modes $\alpha=(\alpha_0,\alpha_1)$ 
reflect oscillations in the variable 
$(\theta_0,\theta_1)=(k x/\ep, - k t/\ep)$.

This choice of Ansatz is motivated by the following.
The oscillations in $k x/\ep$
are anyhow present in the initial data,
and the hyperbolic feature of the equations ensures they are propagated
into oscillations in the variable $(k x/\ep, - k t/\ep)$. Nonlinear interaction of waves next
makes sure that all harmonics are created along time evolution in \eqref{MBeps}. 
To be more precise, we show that all oscillations that are characteristic  
for the Maxwell-Bloch system are propagated.
They are given, {\em for 
fields and populations}, by the characteristic set
\begin{equation*} 
\mathcal{C} := \cp\cup\cm\cup\cz, 
\end{equation*}
where we define 
\begin{equation} 
\label{cpcz} 
\ch_\pm = \left\{ \alpha=(\alpha_0,\alpha_1)\in\Z^{2d}\nz \mid
\alpha_1=\pm \alpha_0 \right\} \quad \text{and} \quad 
\cz = \left\{ \alpha=(\alpha_0,\alpha_1)\in\Z^{2d}\nz \mid \alpha_1=0\right\},
\end{equation}
corresponding respectively to the waves propagating to the left 
and to the right, and to purely spatial oscillations. 
The reader should be cautious about the fact that the Fourier mode $\alpha=0$ is {\em not} considered as a part
of this characteristic set.
These sets essentially correspond to the characteristic variety of the linear Maxwell part of the equations.
{\em For coherences}, we show the characteristic frequencies
are those which resonate with 
some transition energy $\omega(m,n)$. They are given,
for any wavenumber $k \in \Z^d$, by the resonant set
\begin{align}
\label{reso}
\res = \{ (m,n,\alpha_0,\alpha_1) \in \{1,\dots,N\}^2\times\Z^{2d}
\mid \vf\cdot\alpha_1= \omega(m,n) \}.
\end{align}

\paragraph{Describing the asymptotic dynamics  (2) - Rectification effects.}
\mbox{} The \-rec\-tification phenomenon is 
the creation of non-oscillating terms through the nonlinear interaction of oscillating terms. 
This phenomenon is typical of quadratic systems such as \eqref{MBeps}.

It is the reason for our introduction of the intermediate time scale $T=t/\sqrt{\ep}$ in (\ref{momo1}).
This scale  captures
the evolution of the system between the 
macroscopic length
$|x|=\mathcal{O}(1)$ and the wavelength~$\mathcal{O}(\ep)$.
The first need for the intermediate scale $T$ is, the non-oscillating terms
induce a 
secular growth of the necessary corrector terms (see 
\cite{Lan03}). In that circumstance, it turns out that imposing a \emph{sublinearity condition},
\begin{equation} \label{sslin}
\frac{1}{T} \| \sol^j \|_{L^\infty_{t,x,y,z,\sigma,\theta}} 
\tendlorsque{T}{+\infty} 0 \mbox{ for } j=1,2,
\end{equation}
ensures smallness of correctors in \eqref{sola}, namely
$\rep \, \sol^j(t,x,y,z,t/\rep,\sigma,\theta)=o(1)$ when 
$\ep\rightarrow0$ as desired ($j=1,2$). A second reason for the introduction of the scale $T$ is that,
already in the less singular regime of 
diffractive optics, the leading profile $\sol^0$ obeys 
some linear, constant coefficients hyperbolic evolution system 
with respect to time $T$, and the $T$-sublinearity condition for correctors 
provides a unique way of determining  $\sol^0$ 
(see \cite{JMR98}, \cite{Lan98}, and \cite{Dum03}, 
\cite{Dum04} for the variable coefficients case), \via a detailed analysis of wave interactions at the  
scale $T$ (see Section~\ref{SecInterm}). A similar phenomenon occurs in the present situation as well.

More precisely, the above mentioned system at scale $T$
here takes the following form.
Separating the average part and the oscillating part as
(we use the notation introduced in \eqref{momo1}-\eqref{momo3})
\begin{align*}
{\bf u}^0={\bf u}^0_0+{\bf u}^0_{\rm osc},
\qquad
(\text{which defines  } \fzo\equiv\fiel^0-\fiel^0_0),
\end{align*}
we prove that ${\bf u}^0_0$ satisfies a system of the form
\begin{equation} 
\label{interm1} 
\dT\fiel^0_0+M_2(0,\partial_y,\partial_z)\fiel^0_0=0,
\end{equation}
where $M_2(0,\partial_y,\partial_z)$ is a
\emph{matrix}-coefficient differential operator of size $3\times3$. 
Eventually, our analysis shows that \eqref{interm1} has to be solved together 
with the one giving the corrector $\fiel^1_0$, namely
\begin{equation} \label{slow}
\dT \fiel^1_0+M_2(0,\partial_y,\partial_z) \fiel^1_0 = -\mun\fiel^0_0 ,
\end{equation}
where $\mun$ is another matrix differential operator. We establish
that equation (\ref{slow}) possesses a unique solution $ \fiel^1_0$ 
provided $T$-sublinearity of the right-hand-side is imposed. 
It turns out that the same ``secular growth'' analysis is also 
necessary for ${\bf u}^0_{\rm osc}$ and ${\bf N}^0$.

\begin{rem} \label{remcommut}
In the particular case of profiles that {\em do not} depend on $T$, 
the coupled system \eqref{interm1}, \eqref{slow} is actually 
overdetermined (see Remark~\ref{noT}).
Solvability is only recovered provided the matrix operators $\mun$ 
and $M_2(0,\partial_y,\partial_z)$ commute. This very particular 
and important situation occurs in the Transverse Magnetic case 
(see below), where these operators are all scalar. In the 
general case, rectification enforces the introduction of time 
$T$ to make the set of profile equations solvable.
\end{rem}

\paragraph{Describing the asymptotic dynamics  (3) - the coherences.}
For coherences, we establish the relation 
$$
\coh^0\equiv{\bf C}^{0,1} \, e^{-\sigma}=
\left(\sum_{\alpha\in\Z^{2d}} {\bf C}^{0,1}_\alpha e^{i \alpha \cdot \theta}\right) \, 
e^{-\sigma}.
$$
In other words, at dominant order coherences decay as $\exp(-\gamma t/\ep)$. Besides,
the coefficients ${\bf C}^{0,1}_\alpha$ are shown to satisfy (see \eqref{intermtranspCoh0}) 
\begin{align}
\label{coco1}
\forall(m,n,\alpha)\in\res, \qquad 
\dT \coh^{0,1}_{m,n,\alpha} = 
i\,[E^{0}\cdot\Gamma , \coh^{0,1}]_{m,n,\alpha}.
\end{align}
The set $\res$ is defined in \eqref{reso}.
The other components $\coh^{0,1}_{m,n,\alpha}$, for which $(m,n,\alpha)\notin\res$, 
are shown to vanish (and $E^{0}$ is determined independently: see the 
paragraph ``Field dynamics'' below).

Now the resolution of (\ref{coco1}) cannot rely on the above mentioned sublinearity condition,
due to the fact that a generic solution to (\ref{coco1}) grows like $\exp(K T)$ or so at least, for some constant $K>0$.
Similarly, the analogous equations on the correctors
$\coh^1$ and $\coh^2$ cannot be solved using the sublinearity condition neither.
As a consequence, the evolution of $\coh^0$ remains
a priori undetermined.
The key point now comes from the off-diagonal relaxations. Indeed,
we do prove that $\czu$ satisfies a bound of the form
$$
 \forall\mu\in\N^{5+2d},  \quad
\left|\partial_{t,T,x,y,z,\theta}^\mu\czu(T)\right|
\leq K_1e^{K_2T},
$$
for some $K_1, K_2>0$, so that the product
$
\left(\czu \, \exp(-\sigma)\right)\big|_{T=t/\rep, \sigma=-\gamma t/\ep}
$
has size
$
\mathcal{O}\Big(\exp(K_2 t/\rep) \times \\
\times \exp(-\gamma t/\ep)\Big),
$
hence is negligible when $\ep\to 0$, as desired. 
We notice \emph{a posteriori} that the lack of knowledge in 
the evolution of $\coh^0$ is harmless, since 
coherences live during an initial layer of size $\mathcal{O}(\ep)$ 
only. 
We also show that the correctors $\coh^1$ and $\coh^2$ 
are negligible thanks to estimates in the same vein 
(see Section~\ref{SecExp}).

\paragraph{Describing the asymptotic dynamics  (4) - the electromagnetic field.}
First, we show that ${\bf u}^0\equiv{\bf u}^{0,0}$, so that the 
electromagnetic field  does not undergo the same decay as coherences. 
Next, we prove that $\fiel^0$, conveniently decomposed into its average 
and its oscillatory part, may be written 
\begin{equation}
\label{decfi} 
\fiel^0=
\underbrace{\fiel_{0,0}^0+\fiel_{0,+}^0+\fiel_{0,-}^0}_{=\fiel^0_0}
+
\underbrace{\fiel^0_{\rm time}+ \fiel^0_{\rm space}}_{=\fiel^0_{\rm osc}},
\end{equation}
where we define
the purely spatial oscillations,
and temporal oscillations\footnote{to be accurate, these oscillations, 
namely $\fiel^0_{\rm time}$, involve both  time {\em and} 
space variables - we neertheless keep the denomination "time" 
for this part of the oscillations.}, as
\begin{align*}
\fiel^0_{\rm space}
=
\sum_{\alpha\in\cz} \fiel^0_\alpha \, e^{i\,\alpha\cdot\theta},
\qquad
\fiel^0_{\rm time}
= 
\sum_{\alpha\in\cp\cup\cm} \fiel^0_\alpha \, e^{i\,\alpha\cdot\theta}.
\end{align*}
The three mean terms $\fiel_{0,0}^0$, $\fiel_{0,+}^0$ and $\fiel_{0,-}^0$ are defined in the course of the analysis.
The above decomposition entails the fact that we have 
${\bf u}^0_\alpha=0$ whenever $\alpha\notin\cz\cup\cp\cup\cm$. 
This is a polarization condition, as usual in optics. Besides, there are algebraic projectors $\pi_0$, $\pi_+$ and 
$\pi_-$, such that for any $\alpha\in\ch_\iota$, we have $\pi_\iota\fiel^0_\alpha=\fiel^0_\alpha$ ($\iota=0,+,-$).

We also show that functions $\fiel^0_{0,0}$,  
$\fiel^0_{\rm space}$, $\fiel^0_{0,+}$, $\fiel^0_{0,-}$ 
do not depend on the slow time $t$, {\em i.e.}
\begin{equation}
\label{bahbouh}
\dt\fiel_{0,0}^0=0, \quad \fiel^0_{\rm space}=0, \quad 
\dt \fiel_{0,\pm}^0=0,
\end{equation} 
and they satisfy, from \eqref{interm1},
\begin{equation}
\label{nltran}
\dT\fiel_{0,0}^0=0, \quad \dT \fiel^0_{\rm space}=0, \quad 
\left( \dT\pm\sqrt{\Dy^2+\Dz^2} \right) \fiel_{0,\pm}^0=0. 
\end{equation}
The latter is a {\em two-dimensional} (non local) transport 
equation at time scale $T$. The time oscillations 
$\fiel^0_{\rm time}$ satisfy
\begin{align} 
\label{22}
& \forall\alpha\in\ch_\pm, \quad \dT 
\fiel^0_\alpha = 0 , \quad \text{and}
\\
\label{23}
& \qquad (\dt  
\pm \dx)\,\fiel^0_\alpha \mp \frac{i}{2\vo\cdot\alpha_0} 
\Delta_{y,z} \fiel^0_\alpha =
\pi_\pm \Big( 0 , i\,\Tr 
  \left( \Gamma (\Og (\Og-\vf\cdot\alpha_1)^{-1} -1) 
  [E^0_{\rm time}\cdot\Gamma,\pop^0]_\alpha \right) \Big),
\end{align}
a nonlinear 
Schr\"odinger equation with respect to time $t$, describing diffraction in the transverse
variable $(y,z)$, and quadratic coupling between light and matter.
Technically, a standard Diophantine hypothesis ensures that 
$\vo\cdot\alpha_0$ is ``not to small'' in \eqref{23} -- see 
Hypothesis~\ref{HypDioph}.

\paragraph{Describing the asymptotic dynamics  (5) - the populations.}
As for the electromagnetic field, the dominant term 
${\bf N}^0\equiv {\bf N}^{0,0}$ in the populations does not decay 
in time at first order. In addition, we show that ${\bf N}^0$ 
only possesses spatial oscillations, thanks to the identity
$$
\pop^0=\pop^0_0+\pop^0_{\rm space}.
$$
This is a polarization property. Furthermore, $\pop^0$ 
does not depend on the intermediate time $T$,
\begin{equation}
\label{dtnz}
\dT \pop^0=0.
\end{equation}
Lastly, at the macroscopic time scale $t$,
populations evolve according to a Boltzmann-type 
equation, with transition rates that are the sum of the linear rates $W$ (see (\ref{wdiese})), and of
nonlinear rates that depend quadratically on $E^0$ and change with the  
frequency $\theta$. Precisely, we get the relation 
\begin{align}
\nonumber
\dt  \pop^0 = W\op \pop^0 
& - \Big[ (E^0_{\rm osc} + E^0_{0,0})\cdot\Gamma ,
(i\,\Og-\vf\cdot\dthu)^{-1} \crod{(E^0_{\rm osc} +
E^0_{0,0})\cdot\Gamma}{\pop^0} \Big]_{\rm d,space}
\\
\label{rate} 
& - \Big[ (E^0_{\rm osc} + E^0_{0,0})\cdot\Gamma ,
(i\,\Og-\vf\cdot\dthu)^{-1} \crod{(E^0_{\rm osc} +
E^0_{0,0})\cdot\Gamma}{\pop^0} \Big]_{\rm d,0} .
\end{align}
We refer to \eqref{slowtranspPop0} below. In practice, this equation is of the form
\begin{align*}
\dt 
& \pop^0_\alpha = W\op \pop^0_\alpha +\sum_{\alpha'+\alpha''=\alpha}
\sum_\beta W(\alpha',\alpha'',E^0_{\alpha'},
E^0_{\alpha''-\beta})\op \, \, \pop^0_{\beta},
\end{align*}
for some Pauli coefficients $W(\alpha',\alpha'',E^0_{\alpha'},
E^0_{\alpha''-\beta})$ that depend on the Fourier indices
$\alpha'\in\Z^{2d}$,
$\alpha''\in\Z^{2d}$,
$\beta\in\Z^{2d}$, as well as linearly on each variable
$E^0_{\alpha'}$ and $E^0_{\alpha''-\beta}$. 

\paragraph{An illustrative example : the Transverse Magnetic 
(TM) case.}

We detail here the above analysis in the classical TM case 
(see \cite{NM92}), for which formulae become more explicit, 
and the geometry of the problem is a bit simpler. The TM case 
is a 2-dimensional case, for which in addition fields have 
the particular polarization
$$
B=B(t,x,y)=
\begin{pmatrix}
B_x \\
B_y \\
0
\end{pmatrix},
\quad
E=E(t,x,y)=
\begin{pmatrix}
0 \\
0 \\
E
\end{pmatrix}.
$$
Firstly, we give the explicit profile equations in this 
specific situation. Secondly, we show that the introduction of 
the intermediate time $T$ is not necessary in that case 
(see Remark~\ref{remcommut}). 
Finally, we perform the analysis for \emph{prepared data}, 
{\em i.e.} when coherences vanish at leading order. 
Eventually, we show that the approximation is better than 
in \eqref{stab}, in that {\bf Theorem~\ref{ThStabTM}} asserts
$$
\forall\mu\in\N^2, \quad 
\|\partial^\mu_{x,y}\left(\sole-\sola\right)\|
_{L^\infty([0,t_\star]\times\R^2)}=\mathcal{O}(\rep) .
$$

\subsection{Outline of the paper}

In Section~\ref{SecAnsatz}, we describe the Ansatz 
\eqref{sola}, the choice of scales and phases. We also discuss the needed 
Diophantine hypothesis~\ref{HypDioph}  on the wave numbers $(k^1,\ldots,k^d) \in \R^d$ involved in
(\ref{wn}). In Section
~\ref{SecFormal}, we proceed with the construction of an 
approximate solution ${\bf U}^\ep_{\rm app}={\bf U}^\ep_{\rm app}(t,x,y,z,T,\theta)$
that is  consistant with Maxwell-Bloch's system 
\eqref{MBeps}. Sections~\ref{SecAlg} and \ref{SecFast} 
contain the fast scale analysis, {\em i.e.} the analysis in the variable $\theta$. 
This is a Fourier analysis, leading to the usual 
characteristic sets and group velocity of geometric optics. 
Section~\ref{SecInterm} is devoted to the intermediate 
scale analysis of the resulting 
profile equations at the time scale $T$. The sublinearity condition \eqref{sslin} 
is involved (in Section~\ref{SecAvOp}), that allows to treat wave 
interactions. We also perform (in Section~\ref{SecExp}) 
the \emph{ad hoc} analysis for the coherences ${\bf C}^\ep_{\rm app}$, which have 
exponential growth in $T$ and violate the sublinearity condition. In Section
~\ref{SecExist}, we eventually solve the Cauchy problem for profile 
equations in all variables $t$, $T$,  and $\theta$. To complete the analysis, Section~\ref{SecCV} provides the proof of 
our main result, namely that ${\bf U}^\ep_{\rm app}$ indeed  approximates the true solution
${\bf U}^\ep$ to the Maxwell-Bloch system \eqref{MBeps}, see equation \eqref{stab} above.
The precise form of our result is given in Theorem~\ref{ThStab}. Lastly, Section
~\ref{SecTM}  goes again through the whole analysis, yet in the 
simpler Transverse Magnetic case, and for prepared initial 
data. The corresponding result is given in Theorem~\ref{ThStabTM}.

\medskip

Our main theorems are 
 Theorem~\ref{ThStab} and
 Theorem~\ref{ThStabTM}.

\section{Formulating the Ansatz} \label{SecAnsatz}

We solve the Cauchy problem associated with \eqref{MBeps} 
(with unknown $\sole=(\fiele,\cohe,\pope)$, $\fiele=(B^\ep,E^\ep)$) 
for rapidly oscillating initial data of the form
\begin{equation} \label{inidata}
\sole_{\rm ini}(x,y,z) =
\sum_{\beta\in\Z^d} \phi^\ep_\beta(x,y,z) \,
\exp\left(i \, (\beta\cdot\vo) \, \frac{x}{\ep}\right),
\end{equation}
where the wave vector $k$ is
$
\vo=(\vo^1,\ldots,\vo^d)\in\R^d$
for some 
$
d\in\N^\star,
$
its coordinates are chosen 
$\Q$-independent, and the scalar product 
$\beta \cdot k$ in \eqref{inidata} denotes
$
\beta \cdot k = \beta^1 k^1 + \cdots + \beta^d k^d.
$
The $\Q$-independence of $k$'s coordinates ensures that for
$\beta\in\Z^d$, relation $\beta \cdot k =0$ holds if and only 
if $\beta =0$. 
The vector $k$ collects the {independent} oscillations (in $x$) 
carried by the initial signal 
$\left( B^\ep,E^\ep,\rho^\ep \right)_{\rm ini}(x,y,z)$, while
the integer $\beta^j \in \Z$
numbers the various {harmonics} corresponding to the 
phase $k^j \,x /\ep$ ($j=1,\ldots,d$). 
We choose an initial signal that carries all possible harmonics.
This harmless choice is motivated by the fact that the original, nonlinear,
Maxwell-Bloch system 
anyhow generates nonlinear wave interaction, which implies that any
initial oscillation creates the associated harmonics along wave propagation.

The off-diagonal relaxation term $- \gamma {\bf C}^\ep/\ep$ 
in \eqref{MBeps} enforces exponential decay of
coherences $\cohe$, so that, at first order, 
we expect only the fields $\fiele=(B^\ep,E^\ep)$ and
the populations $\pope$ to propagate. Hence space and time 
oscillations are expected to
be generated from the initial spatial ones only through 
Maxwell's equations and
through the equation for the populations. 
In a similar spirit, since our scaling 
postulates that variations in the transverse variables $(y,z)$ 
only occur at
the scale $\rep$ (and not $\ep$), we also expect that oscillations will not occur in the $(y,z)$ variables. In other
words, the relevant Maxwell equations for
oscillations are expected to be 1D in the $x$ direction.
The corresponding
characteristic variety, in Fourier variables, reads
$$
\mathcal{C}_{Maxwell} \cup \mathcal{C}_{populations}
=
\left\{
(\tau,\xi) \in \R^2\nz \mid \tau=0 \mbox{ or }
\tau^2=\xi^2
\right\}.
$$
Eventually, propagation of space and time
oscillations are expected to occur \emph{via}
the collection of phases 
$$
k^j \, x - \omega^j \, t, 
\quad\mbox{ with } \, \, 
\omega^j=0 \,  \mbox{ or } \, \omega^j= \pm \,k^j
\quad ( j=1,\ldots, d).
$$ 
More precisely, the relevant oscillations in our analysis 
are expected to be all the harmonics
$\beta \cdot (k \, x - \omega \, t)$, 
as the multi-index $\beta$ runs in $\Z^d$, and the $j$-th coordinate 
of $\omega$ is either $0$, $k^j$ or $-k^j$. 
An instant of reflexion shows that this ensemble  coincides
with the collection of phases 
$(\alpha_1 \cdot k) \, x - (\alpha_0 \cdot k) t$ as $\alpha_0$ 
and $\alpha_1$ run in $\Z^d$. This motivates the following
\begin{nota}
For any multi-index $\alpha=(\alpha^1,\ldots,\alpha^{2d})\in \Z^{2d}$,
we decompose
\begin{align*}
\alpha:=(\alpha_0,\alpha_1), \text{ with } 
& \alpha_0=(\alpha^1,\ldots,\alpha^d)\in\Z^d,
\text{ and } 
\alpha_1=(\alpha^{d+1},\ldots,\alpha^{2d})\in\Z^{d}.
\end{align*}
Accordingly, any (smooth enough) function ${\bf V}(\theta)$ defined 
over the  torus $\T^{2d}$ may be decomposed in Fourier
series as
\begin{align*}
{\bf V}(\theta)
=
\sum_{\alpha\in\Z^{2d}}
{\bf V}_\alpha \, \exp\left(i \, \alpha \cdot \theta \right)
=
\sum_{\alpha\in\Z^{2d}}
{\bf V}_\alpha \,
\exp\left(i \, (\alpha_0\cdot\theta_0 + \alpha_1\cdot\theta_1) \right).
\end{align*}
\end{nota}
With this notation, \emph{resonances} are identified by the
characteristic set $\mathcal{C} = \cz\cup\cp\cup\cm$, 
where $\cz$, $\cp$, and $\cm$ have been defined in \eqref{cpcz}. 

\medskip

Due to the fast decay of coherences, another set of phases needs 
to be introduced.
Namely, the off-diagonal relaxation term $- \gamma {\bf C}^\ep/\ep$ 
in \eqref{MBeps} leads to introduce the complex phases 
$\exp(-\gamma t/\ep)$ as well as all its harmonics
$\exp(-\kappa \gamma t/\ep)$, where $\kappa\in\N$.

\medskip
With all these considerations in mind,  we can now define
the approximate solution we seek, as
\begin{align}
\label{solapp}
\nonumber
&
\quad
\sola = (\fiela,\coh\ea,\pop\ea) ,
\quad \fiela = (B\ea,E\ea),
\quad
\text{where}
\\  
&
\quad \sola(t,x,y,z) =
\sum_{j=0}^2 \rep^j\sol^j(t,x,y,z,T,\sigma,\theta_0,\theta_1)\Bigg|_{
T={t}/{\rep},
\sigma={\gamma t}/{\ep},(\theta_0,\theta_1)=\left(\vo x , -\vf t\right)/\ep},
\\
\nonumber
&
\text{and for any $j=0,1,2$, and any } \sigma\geq0,T\geq0,~\theta\in\T^{2d}, \quad \text{we set}
\\ 
\nonumber
&
\quad
\sol^j\varTdt
=
\sum_{\alpha\in\Z^{2d}} \sum_{\kappa\in\N}
\sol^{j,\kappa}_\alpha\varT \,
\exp\left( i \alpha \cdot \theta \right) \exp(-\kappa\sigma).
\end{align}

\begin{rem} \label{RemAnsatz}{\bf (discussion of the chosen profiles).}\\
(i) Coming back to the original scales of the problem, we may
represent the solution $\sole$ under the form 
$$\sole\var=\tilde{\sol}^\ep(t,x,\rep y,\rep z),$$
where the function $\tilde{\sol}^\ep$ has variations in
$\tilde{y}=\rep y$ and $\tilde{z}=\rep z$ at scale $\rep$, and solves
\eqref{MBeps} with the operator $\roteps$ replaced with $\curl$, while the associated 
initial datum is of the form
$$\tilde{\sol}^\ep_{|_{t=0}}(x,\tilde{y},\tilde{z})=\sole_{|_{t=0}}
\left(x,\frac{\tilde{y}}{\rep},\frac{\tilde{z}}{\rep}\right)=:
\solc_{|_{t=0}}\left(x,\frac{\tilde{y}}{\rep},\frac{\tilde{z}}{\rep},\frac{x}{\ep}\right).$$
This is the 3-scales setting of \cite{Dum04} (where only the case of
quasilinear, non-dispersive systems is adressed). It leads to a
profile representation of the form
$$\tilde{\sol}^\ep(t,x,\tilde{y},\tilde{z})= 
\tilde{\solc} \left( t,x,\tilde{y},\tilde{z},\frac{\psi}{\rep},
\frac{\phi}{\ep} \right).$$  
where $\phi=\phi(t,x,\tilde{y},\tilde{z})$ is the collection of the two phases 
$kx$ and $-kt$, while  $\psi=\psi(t,x,\tilde{y},\tilde{z})$ is
a collection of ``intermediate phases'' (in \eqref{solapp},
$\tilde{\solc}$ does not depend explicitly on $\tilde{y},\tilde{z}$,
so that $\psi=(t,\tilde{y},\tilde{z})$). When non-oscillating terms
are present at first order ($\tilde{\sol}_0^0 \neq 0$), it is in
general necessary to put in this collection $\psi$ an intermediate
time $t$ in order to solve the profile equations \emph{via} the
analysis of Section~\ref{SecInterm}. This is the reason
why the variable $T=t/\rep$ is present in \eqref{solapp} (it may be
unnecessary: see \cite{Dum04}, Remark~1.4, and the Transverse
Magnetic case below, Section~\ref{SecTM}, for which
$\psi=\tilde{y}$). This intermediate time $T$ captures the evolution
of intermediate scales  $\rep$, between ``macroscopic length''
$\mathcal{O}(1)$ and wavelength $\ep$.
\\
(ii) The consistancy of the chosen Ansatz requires that $\rep\sol^1 \ll \sol^0$ whenever $\sigma=\gamma
t/\ep$, $T=t/\rep$,  $\theta=(\vo x/\ep , -\vf t/\ep)$. This requirement
enforces $T$-sublinearity of the non-exponentially decaying part (corresponding to $\kappa=0$) of 
the correctors $\sol^{j,\kappa}$ ($j\geq1$), see \eqref{sublin}.
The  exponentially
decaying correctors ($\kappa\geq1$) may be of the same order as the first 
profile $\sol^{0,\kappa}$: they anyhow lead to $o(\rep)$ error terms, see 
Section~\ref{SecExp}. 
\\
(iii) The analogy with \cite{Dum04} shows that we could treat the
same problem with {\em curved} ({\rm i.e.} nonlinear) phases $\phi$.
This situation arises in particular when dealing with  {\em inhomogeneous media}, where variable
coefficients ({\rm i.e.} variable electric and magnetic permittivity) are  involved in the  Maxwell-Bloch system. 
\\
(iv) We emphasize the fact that in the present work, we deal with
\emph{large population} variations, {\em of order one}.

In \cite{JMR00}
(long time diffraction, p. 248), with no relaxation terms,
{\em transparency} allows to transform  the original Maxwell-Bloch system, \emph{via} 
a change of dependent variables, into a dispersive and quadratic system 
of the form $L(\ep\partial) U = \ep^2 f(U)$, with
$U=(U^I,U^{I\!I})$, where $(B,E,\coh)=\ep U^I$, and
$\pop=\ep^2~U^{I\!I}$. The system is considered over times of order
$1/\ep$. Taking relaxations into account, we consider here a similar system, at the same
space and time scales (oscillations have frequency $1/\ep$, and
propagation is considered over times of order $1/\ep$), yet in a situation where $(B,E,\pop,\coh)$ are larger,
of order 1, and the quadratic coupling term $f(U)$ is stronger, of size  $\ep$, namely
\begin{align*}
\left\{
\begin{array}{l}
\vspace{0.2cm}
\ep\dt\cohe = -i \crod{\Omega}{\cohe} 
+ i \ep \crod{E^\ep\cdot\Gamma}{\cohe} 
+ i \ep \crod{E^\ep\cdot\Gamma}{\pope} - \gamma \cohe,
\\
\ep\dt\pope = i\ep\crd{E^\ep\cdot\Gamma}{\cohe} + \ep^2W\op\pope.
\end{array} 
\right.
\end{align*}
Of course, the stronger quadratic interaction term is  balanced by the off-diagonal relaxations: note however that relaxation
only affects part of the unknowns (the coherences), and one key aspect of our analysis precisely relies in the analysis
of  the stronger interaction term.

Such large population variations are also considered in \cite{BBCNT04},
without relaxation, yet for waves that are polarized in a specific way,
and in the case of a three-level Bloch system only ({\em i.e.} $N=3$
in our notation). In that case, additional conservation properties 
are at hand. 
On top of that, \cite{BBCNT04} considers a 
weaker coupling term of size $\ep^2$. 
The limiting system obtained in  \cite{BBCNT04} 
is of Schr\"odinger-Bloch type.
\end{rem}

Since small divisor estimates naturally
enter the analysis below, we readily formulate the usual
Diophantine assumption we shall need on the wave-vector $k$.
It will be used in order to invert the differential
operators acting in the $\theta$ variable on the various (smooth)
profiles ${\bf U}^{j,\kappa}\varTt$. We refer for instance to
\cite{JMR93}.
\begin{hyp} \label{HypDioph}
The wave vector  $k\in\R^s$ from \eqref{inidata} is 
{\em Diophantine}, namely
$$
\exists C,a>0, \quad \forall \beta\in\Z^d\nz,
\quad |\beta\cdot k|
\geq C|\beta|^{-a}.
$$
\end{hyp}

\begin{rem}
The above assumption is harmless. Indeed, the following fact is 
well-known. Pick {\em any} exponent $a>d-1$. Then, the set
$$
\{k \in \R^{d} \mid \exists C>0, \quad \forall \beta\in\Z^d\nz,
\quad |\beta\cdot k|
\geq C|\beta|^{-a}\}
$$
has full measure in $\R^d$. 
In other words, almost any $k \in \R^d$ (for the Lebesgue measure)
has the Diophantine property.
\end{rem}

\section{Formal expansions and approximate solution}
\label{SecFormal}

\begin{nota} \label{NotSymbMax}
Denote by $M(\dt,\dx,\dy,\dz)$ the order 1 differential operator in
Maxwell's equations, 
$$
M(\dt,\dx,\dy,\dz) = 
\begin{pmatrix}
0 & \curl \\
-\curl & 0 
\end{pmatrix} =
\dt + A_x\dx + A_y\dy + A_z\dz ,
$$
where the $A_j$'s are $6\times6$ real symmetric matrices. Set also 
\begin{equation*}
\mun=M(\dt,\dx,0,0)=\dt + A_x\dx, \quad 
\md=M(\dT,0,\dy,\dz)=\dT  + A_y\dy + A_z\dz.
\end{equation*}
\end{nota}

With this notation at hand, plugging the Ansatz \eqref{solapp} into the Maxwell-Bloch
system~\eqref{MBeps}, we get that ${\bf U}^\ep_{\rm app}$ satisfies the original Maxwell-Bloch equations
up to  a residual $r^\ep$ that is given as follows (here we use the notation 
\eqref{Og} for $\Og$)

\begin{prop} \label{PropDA}
The residual $$r^\ep := \Lep\sola - F^\ep(\sola)$$ has the profile 
representation
\begin{align*}
& r^\ep\var = \mathcal{R}^\ep(t,x,y,z,T,\sigma,\theta)\Bigg|_{
T={t}/{\rep},
\sigma={\gamma t}/{\ep},(\theta_0,\theta_1)=\left(\vo x , -\vf t\right)/\ep},
\end{align*}
where
\begin{align*}
& \mathcal{R}^\ep = \sum_{j=-2}^3 \sum_{\kappa\geq0}\rep^j 
r^{j,\kappa}\varTt\exp(-\kappa\sigma),
\end{align*}
and the first terms ($-2\leq j\leq0$) are (see Section~\ref{Secres} 
for the others)
\begin{align*}
r^{-2,\kappa} = \left(
\begin{array}{l}
\vspace{0.15cm}
 \mudth\fiel^{0,\kappa} \\
\vspace{0.15cm}
  (i\,\Og-\gamma\kappa-\vf\cdot\dthu)\coh^{0,\kappa} \\
\vspace{0.15cm}
 (-\gamma\kappa-\vf\cdot\dthu)\pop^{0,\kappa} 
\end{array}
\right),
\end{align*}
\begin{align*}
r^{-1,\kappa} = \left(
\begin{array}{l}
\vspace{0.15cm}
 \mudth\fiel^{1,\kappa} + \md\fiel^{0,\kappa} 
- \left( 0 , i\, \Tr(\Gamma\Og\coh^{0,\kappa}) \right) \\
\vspace{0.15cm}
(i\,\Og-\gamma\kappa-\vf\cdot\dthu)\coh^{1,\kappa} 
+ \dT\coh^{0,\kappa} 
-i\,\crod{E^0\cdot\Gamma}{\coh^0+\pop^0}^\kappa \\ 
\vspace{0.15cm}
 (-\gamma\kappa-\vf\cdot\dthu)\pop^{1,\kappa} + \dT \pop^{0,\kappa} 
- i\,\crd{E^0\cdot\Gamma}{\coh^0}^\kappa
\end{array}
\right),
\end{align*}
\vspace{0.4cm}
\begin{align*}
r^{0,\kappa} = \left(
\begin{array}{l} 
 \mudth\fiel^{2,\kappa} + \md\fiel^{1,\kappa} \\
\vspace{0.15cm}
\qquad
+ \mun \fiel^{0,\kappa} 
+ \Big( 0 , -i\, \Tr(\Gamma\Og\coh^{1,\kappa}) 
+ i\, \Tr(\Gamma [E^0\cdot\Gamma,\coh^0+\pop^0]^\kappa) \Big) \\
 (i\,\Og-\gamma\kappa-\vf\cdot\dthu)\coh^{2,\kappa} + \dT \coh^{1,\kappa} 
+ \dt \coh^{0,\kappa} 
\\
\vspace{0.15cm}
\qquad
- i\,\crod{E^0\cdot\Gamma}{\coh^1+\pop^1}^\kappa 
- i\,\crod{E^1\cdot\Gamma}{\coh^0+\pop^0}^\kappa  \\
\vspace{0.15cm}
  (-\gamma\kappa-\vf\cdot\dthu)\pop^{2,\kappa} + \dT \pop^{1,\kappa} 
+ \dt \pop^{0,\kappa} - i\,\crd{E^0\cdot\Gamma}{\coh^1}^\kappa 
- i\,\crd{E^1\cdot\Gamma}{\coh^0}^\kappa - W\op \pop^{0,\kappa}
\end{array}
\right).
\end{align*}
\end{prop}

With this computation at hand, we wish to construct the first profile
$\sol^{0,\kappa}=(\fiel^{0,\kappa},\coh^{0,\kappa},\pop^{0,\kappa})$,
in $\mathcal{C}^0([0,t_\star]_t,
\mathcal{C}^0([0,+\infty[_T,H^\infty(\R^3\times\T^{2d})))$ 
for some $t_\star>0$, together with \emph{correctors}
$(\fiel^j,\coh^j,\pop^j)$ ($j=1,2$), in such a way that 
the residual $r^\ep$ is small (\cf Section~\ref{SecStab}). More precisely,
we shall impose $r^{-2}=r^{-1}=r^0=0$,
and show that this procedure completely
determines $\sol^0$. This is obtained  by
decomposing
successively (with increasing $j$) the equations obtained for the
$\sol^j$'s, and separating characteristic and noncharacteristic modes
$(\alpha,\kappa)$. The necessary linear algebra tools are developed in the next paragraph.

\subsection{Rapid modes and algebraic projectors} \label{SecAlg}

\paragraph{$\bullet$ Tools needed to deal with the electromagnetic field.}
The Fourier series representation, in the $\theta$ variable, of the first 
equation stemming from $[r^{-2,\kappa}=0]$ is (according to 
notation~\ref{NotSymbMax})
$$
\muda\fiel^{0,\kappa}_\alpha=0,\quad\forall\alpha\in\Z^{2d}.
$$
For each $\alpha$, this is a system of linear equations in 
$\C^6$, with matrix 
$$\muda =  i\vf\cdot\alpha_1 + i(\vo\cdot\alpha_0)A_x,$$
which may be singular only if $\kappa=0$.
To deal with this equation,  we need the following classical definitions and lemmas 
(\cite{Lax57}, \cite{JMR93}, \cite{DJMR95}, \cite{JMR98}).

\begin{defn} \label{DefProj}
For each $\alpha\in\Z^{2d}$, let $\pi_\alpha$ be the 
orthogonal projection in $\R^6$ onto the kernel of $M_1(-i k \cdot \alpha_1,i k \cdot \alpha_0)$,. 
Denote by $M_1(-i k \cdot \alpha_1,i k \cdot \alpha_0)^{-1}$ the inverse of $M_1(-i k \cdot \alpha_1,i k \cdot \alpha_0)$, when restricted to 
the space orthogonal to its kernel, namely $M_1(-i k \cdot \alpha_1,i k \cdot \alpha_0)^{-1}$ acts on 
$\ran(1-\pi_\alpha)$. Finally, define the projector $\Pi$ 
on the space of Fourier series, 
$$
\Pi \left(\sum_{\alpha\in\Z^{2d}} u_\alpha e^{i\,\alpha\cdot\theta} \right)
:=\sum_{\alpha\in\Z^{2d}} \pi_\alpha u_\alpha 
e^{i\,\alpha\cdot\theta}.
$$
\end{defn}

A straightforward computation establishes the
\begin{lemme}
For each $\alpha\in\Z^{2d}$, the projector
$\pi_\alpha$ is a homogeneous 
function of $\alpha$ of degree zero. It takes a constant 
(matrix) value $\pi^\iota$ on each component $\ch_\iota$
of the characteristic set ($\iota=+,-,0$), the value one (or identity) 
for $\alpha=0$, and vanishes else. In particular, we have
the identity
\begin{equation}
\label{pipi}
\Pi 
\left(\sum_{\alpha\in\Z^{2d}} u_\alpha e^{i\,\alpha\cdot\theta} \right)
:=
u_0
+
\sum_{\alpha\in\cp} \pi^{+} u_\alpha 
e^{i\,\alpha\cdot\theta}
+
\sum_{\alpha\in\cm} \pi^{-} u_\alpha 
e^{i\,\alpha\cdot\theta}
+
\sum_{\alpha\in\cz} \pi^{0} u_\alpha
e^{i\,\alpha\cdot\theta}.
\end{equation}
In any circumstance, for any Fourier series 
$\displaystyle u = \sum_{\alpha\in\Z^{2d}} u_\alpha 
e^{i\,\alpha\cdot\theta}$, we have
$$\muth u = 0  \Longleftrightarrow \Pi \, u = u.$$
\end{lemme}

Next, the following two lemmas are classical for 
geometric and diffractive optics  with smooth characteristic 
varieties (\cite{Lax57}, \cite{JMR93}, \cite{DJMR95}, and 
also \cite{Texier04} for an elegant unified version). They 
express that the operators acting on the oscillating part 
of the fields are in diagonal form. 
\begin{lemme}[group velocity] \label{LemVg}\mbox{}\\
For any Fourier series 
$\displaystyle u=\sum_{\alpha\in\Z^{2d}\nz} u_\alpha 
e^{i\,\alpha\cdot\theta},$ containing no mean term,
we have, with $\val = \iota$ for $\alpha\in\ch_\iota$ ($\iota=0,+,-$), the relations
\begin{equation*}
\Pi A_y \Pi \, u = \Pi A_z \Pi \, u = 0 , \quad \mbox{ and } \quad
\Pi A_x \Pi \, u = \sum_{\alpha\in\mathcal{C}} \val \pi_\alpha 
u_\alpha e^{i\,\alpha\cdot\theta} =: \vDth \Pi \, u.
\end{equation*} 
\end{lemme}
\begin{lemme}[diffraction]
\label{LemDiff}\mbox{}\\
For any Fourier series 
$\displaystyle u=\sum_{\alpha\in\Z^{2d}\nz} u_\alpha 
e^{i\,\alpha\cdot\theta},$ containing no mean term, we have
\begin{align*} 
&
\Pi \mdyz 
\muthi \mdyz \Pi \, u 
\\
&
\qquad\qquad
= i \,\sum_{\alpha\in\cp\cup\cm}\aal\Delta_{y,z}\pi_\alpha 
u_\alpha e^{i\,\alpha\cdot\theta}
=:i\,\aDth \Delta_{y,z} \Pi 
\, u ,
\end{align*}
where $\displaystyle \aal = \pm \frac{1}{2\vo\cdot\alpha_0}$ 
for $\alpha\in\ch_\pm$ and $\aal=0$ for $\alpha\in\cz$.
\end{lemme}
\mbox{}\\
\noindent{\bf Proof.}
By the definition of the operators involved, for each $\alpha$, one has  
\begin{equation*}
\pi_\alpha \mdyz \muai \mdyz \pi_\alpha 
= \sum_{j,k \in \{y,z\}} \pi_\alpha A_j \muai A_k 
\pi_\alpha \partial_j \partial_k . 
\end{equation*}
Now, $\pi_\alpha$ and $\muai$ are the evaluations, at
$\xi=(\vo\cdot\alpha_0,0,0)$, of the spectral projector
$P(\tau(\xi),\xi)$ resp. of the pseudo-inverse $M^{-1}(\tau(\xi),\xi)$ 
of the complete Maxwell symbol $M(\tau(\xi),\xi)$. Here and in the sequel, we use
the notation
\begin{align*}
\tau_\pm(\xi)=\mp|\xi| \text{ \,  for \, } \alpha\in
\ch_{\iota,\pm}:=\ch_\iota\cap\{\vf\cdot\alpha_1\gtrless0\}, \quad
\iota=\pm,
\end{align*}
together with
\begin{align*}
\tau(\xi)=0 \text{ \,  for \, } \alpha\in\ch_0.
\end{align*}
The quantities $\tau_\pm(\xi)$ and $\tau(\xi)$ are naturally the eigenvalues of the Maxwell symbol $M(\tau(\xi),\xi)$.
On the other hand, we 
have the identity, valid  for any $j,k \in\{x,y,z\}$, see \cite{JMR93} or \cite{Texier04}, 
\begin{align*}
&
\frac{\partial^2\tau}{\partial_{\xi_j}\partial_{\xi_k}}(\xi) 
P(\tau(\xi),\xi)
\\
&
\quad
=
P(\tau(\xi),\xi) \, A_j \, M^{-1}(\tau(\xi),\xi) \, A_k \,
P(\tau(\xi),\xi) 
+ P(\tau(\xi),\xi)\, A_k\, M^{-1}(\tau(\xi),\xi)\, A_j\, 
P(\tau(\xi),\xi).
\end{align*}
Using now the obvious relation
\begin{equation*}
\left(\frac{\partial^2\tau_\pm}{\partial {\xi_j}
\partial {\xi_k}}(\xi)\right)_{i,j} 
= \mp \frac{1}{|\xi|} \left( {\rm id} -
  \frac{\xi\otimes\xi}{|\xi|^2} \right) 
= \mp \frac{1}{|\vo\cdot\alpha_0|} \, {\rm Diag}\,(0,1,1)
\quad \mbox{ for } \xi=(\vo\cdot\alpha_0,0,0),
\end{equation*}
inspection of the five cases $\alpha\in\ch_{\iota,\pm}$ ($\iota=\pm$), and $\alpha\in\cz$, leads 
to  the lemma.
\fin

Lastly, in order to distinguish between propagated and non-propagated 
parts of the profiles, we introduce the following splitting, 
refering to the ``oscillating'' part of the fields, to 
the ``time'' oscillations of the density matrix, or to its ``space'' oscillations.

\begin{defn} \label{DefOsc}
For any Fourier series 
$U=\sum_{\alpha\in\Z^{2d}}U_\alpha e^{i\,\alpha\cdot\theta},$ 
we set
$$
U_{\rm osc} := \sum_{\alpha\neq0} U_\alpha 
e^{i\,\alpha\cdot\theta},
\quad U_{\rm space} := \sum_{\alpha\in\cz} U_\alpha
e^{i\,\alpha\cdot\theta},
\quad U_{\rm time} := \sum_{\alpha\in\cp\cup\cm} 
U_\alpha e^{i\,\alpha\cdot\theta}.
$$
\end{defn}

\paragraph{$\bullet$ Tools needed to deal with the coherences.} 
For each mode $(\alpha,\kappa)$,
the second equation from $[r^{-2,\kappa}_\alpha=0]$ reads 
$$
\forall m,n, \quad
(i(\omega(m,n)-\vf\cdot\alpha_1)+\gamma(1-\kappa))
\coh^{0,\kappa}_{m,n,\alpha} = 0.
$$
Here, only $\kappa=1$ is of interest. This justifies the introduction of the 
resonant set $\res$ we defined in \eqref{reso}.
\begin{rem}
In \cite{BBCNT04}, the wave vector $\vo$ is precisely chosen 
so that $\mathcal{R}(\vo)$ be nonempty.
\end{rem}

\paragraph{$\bullet$ Tools needed to deal with the populations.} 
The populations are scalar variables, and the characteristic frequencies are 
simply the $(\alpha,\kappa)$'s belonging to $\cz\times\{0\}$. No additional tool is needed.

\subsection{Profile equations, fast scale analysis} 
\label{SecFast}

\subsubsection{The residual $r^{-2}$}

According to the notation above, the equation $[r^{-2}=0]$ is equivalent to the
polarization conditions
\begin{align} \label{polarFiel}
&
\fiel^0  =  \fzz,  \qquad\qquad \forall\kappa>0, \quad
\fiel^{0,\kappa}=0, \qquad \text{ and } \qquad 
\fiel^0=\Pi \fiel^{0},
\\
\label{polarCoh}
&
\coh^0 = \czu \, e^{-\sigma}, \qquad \mbox{ and }  \qquad 
\forall (m,n,\alpha)\notin \res, \qquad 
\coh^{0,1}_{m,n,\alpha}=0 , 
\\
\label{polarPop}
&
\pop^0 = \pzz, \qquad\quad\;\; \text{ and } \qquad 
\pop^0= \pop^{0}_0+\pop^{0}_{\rm space}.
\end{align}

In order not to overweight notation, we shall from now on systematically
refer to $\fiel^0$ and $\pop^0$ in the sequel, keeping in mind they do coincide with
$\fiel^{0,0}$ and $\pop^{0,0}$.

\subsubsection{The residual $r^{-1}$}

\paragraph{Equations on the field. } 
We first deal with the mode $\kappa=0$. Separating the mean term 
($\alpha=0$) and the other Fourier modes, we find 
\begin{align}
\label{intermtranspFielMoy0}
&
\partial_T \fiel^{0}_0 +M_2(0,\partial_y,\partial_z) \fiel^0_0= 0,
\\
&
\label{intermtranspFielOsc0}
\dT {\bf u}^0_{\rm osc} = 0.
\end{align}
Obtaining the second equation requires to use \eqref{polarFiel} and Lemma~\ref{LemVg}, which imply
$\Pi M_2(0,\partial_y,\partial_z)  \Pi=0$ and ${\bf u}^0_{\rm osc}=\Pi {\bf u}^0_{\rm osc}$.
The next order profile ${\bf u}^{1,0}$ is then seen to satisfy 
\begin{equation} \label{polarFiel10}
(1-\Pi) \fuzo = - \muthi \mdyz {\bf u}^0_{\rm osc}.
\end{equation}

Secondly, when $\kappa\geq1$, the operator $\mudth$ is invertible, and we get the two values
\begin{align} \label{polarFiel11}
&
\fiel^{1,1} = \mudthi(0,i\Tr(\Gamma\Og\coh^{0,1})),
\\
&
\label{polarFiel11+}
\forall\kappa>1, \quad \fiel^{1,\kappa}=0.
\end{align}

\paragraph{Equations on the coherences.}
We first note that whenever $\kappa\neq1$ the operator $(i\,\Og-\gamma\kappa-\vf\cdot\dthz)$
is invertible. Thanks to
\eqref{polarFiel}, \eqref{polarCoh} and \eqref{polarPop}, this gives the two values
\begin{align} \label{polarCoh10}
&\coh^{1,0} = i\,(i\,\Og-\vf\cdot\dthu)^{-1} 
\crod{E^{0}\cdot\Gamma}{\pop^{0}},
\\
&
\forall\kappa>1, \qquad
\label{polarCoh12+}
\coh^{1,\kappa} = 0. 
\end{align}

When $\kappa=1$, we need to distinguish between resonant and non-resonant 
triples $(m,n,\alpha)$, to obtain
\begin{align} \label{polarCoh11}
&
\forall(m,n,\alpha)\notin\res, \quad 
\coh^{1,1}_{m,n,\alpha} = -(\omega(m,n)-\vf\cdot\alpha_1)^{-1} 
[E^{0}\cdot\Gamma , \coh^{0,1}]_{m,n,\alpha} ,
\\
&
 \label{intermtranspCoh0}
\forall(m,n,\alpha)\in\res, \quad 
\dT \coh^{0,1}_{m,n,\alpha} = 
i\,[E^{0}\cdot\Gamma , \coh^{0,1}]_{m,n,\alpha} .
\end{align}

\paragraph{Equations on the populations.} 
Here we need to distinguish between the values 
 $\kappa=0$ and 
$\kappa\neq0$, as well as between the modes $\alpha\in\cz\cup\{0\}$ and the other Fourier modes.
In that  way we obtain, taking $\kappa=0$ and restricting to modes $\alpha\in\cz\cup\{0\}$, the relation
\begin{equation} \label{intermtranspPop0}
\dT \pop^0 = 0 ,
\end{equation}
together with the following polarization conditions and spectral properties for $\pop^1$ (here we use
(\ref{polarPop}))
\begin{equation} \label{polarPop1}
\forall\alpha\notin\cz\cup\{0\}, \,
\pop^{1,0}_\alpha = 0, \quad
\pop^{1,1} = i\, (-\gamma-\vf\cdot\dthu)^{-1} 
\crd{E^0\cdot\Gamma}{\czu}, \qquad
\forall\kappa>1, \, \pop^{1,\kappa} = 0.
\end{equation}

\subsubsection{The residual $r^0$}

\paragraph{Equations on the field.} 

When $\kappa=0$, we have (using \eqref{polarCoh10} and 
\eqref{polarCoh}) 
\begin{align} 
\label{etoile}
&0=\muth \fdz + \md \fuz + \mun \fzz \\
\nonumber
& \qquad + \Big( 0 , i \, \Tr(\Gamma[E^{0,0}\cdot\Gamma,\pzz]) 
- i \, \Tr(\Gamma i\Og (i\Og-\vf\cdot\dthu)^{-1}
[E^{0,0}\cdot\Gamma,\pzz]_{\rm od}) \Big) .
\end{align}
Again, we separate oscillating and nonoscillating parts in the 
above equation.
In the nonoscillating case, the two nonlinear contributions are 
seen to compensate each other and we obtain
\begin{align} \label{intermtranspFielMoy10}
\md \fiel^{1,0}_0 = - \mun \fiel^{0}_0 , 
\end{align}
On the other hand, the $\Pi$-polarized part of oscillating modes 
turns out to satisfy the relation, using 
\eqref{intermtranspFielOsc0}, \eqref{polarFiel10}, 
\eqref{polarCoh10} and Lemma~\ref{LemVg}, 
and writing $M_2(\dT,\dy,\dz)=\dT+M_2(0,\dy,\dz)$, 
\begin{align*} 
\dT \Pi \fuzo = - \Pi \mun {\bf  u}^0_{\rm osc} 
&
+ \Pi \mdyz \muthi \mdyz  {\bf  u}^0_{\rm osc} \\
&
+ \Pi \Big( 0 , i\,\Tr \left( \Gamma (i\, \Og 
(i\,\Og-\vf\cdot\dthu)^{-1} -1) 
[E^{0}\cdot\Gamma,{\bf N}^0]_{\rm osc} \right) \Big). 
\end{align*}
Thanks to Lemma~\ref{LemVg} and Lemma~\ref{LemDiff}, this is 
rewritten, using ${\bf  u}^0_{\rm osc}=\Pi {\bf  u}^0_{\rm osc}$,
\begin{align} \label{intermtranspFielOsc10}
\nonumber
\dT \Pi \fuzo = 
& - (\dt+\vDth) {\bf  u}^0_{\rm osc} + i\,\aDth \Delta_{y,z} {\bf  u}^0_{\rm osc} \\
&\quad
 + \Pi \Big( 0 , i\,\Tr \left( \Gamma (i\, \Og 
(i\,\Og-\vf\cdot\dthu)^{-1} -1) 
[E^{0}\cdot\Gamma,{\bf N}^0]_{\rm osc} \right) \Big). 
\end{align}
Lastly, applying $(1-\Pi)$ to the oscillating part of equation (\ref{etoile}), we also recover 
(using the fact that 
$\dT \left( (1-\Pi) \fiel^{1,0}_{\rm osc}\right)=0$,
deduced from
equation (\ref{polarFiel10}) together with the identity 
$\dT \fiel^{0}_{\rm osc}=0$),
\begin{align} \label{polarFiel20}
(1-\Pi) \fdzo = 
\nonumber
& - \muthi \Big( \mdyz \fuzo +\mun {\bf  u}^0_{\rm osc} \\
& + \left( 0 , i\,\Tr \left( \Gamma 
    (i\, \Og (i\,\Og-\vf\cdot\dthu)^{-1} -1)
    [E^{0}\cdot\Gamma,{\bf N}^0]_{\rm osc} \right) \right) \Big) .
\end{align}

\vspace{0.2cm}
When $\kappa=1$, using
(\ref{polarFiel11}), (\ref{polarFiel}), (\ref{polarPop}), 
we get in the same way
\begin{align} \label{polarFiel21}
\nonumber
\fdu 
 = - M_1(-\gamma-k\cdot \partial_{\theta_1},k\cdot\partial_{\theta_0})^{-1} 
&
\Big[  \md \muthi 
\left( 0 , i\,\Tr \left( \Gamma (i\, \Og \czu) \right) \right) \\
&
- \left( 0 , i\,\Tr (\Gamma \Og \cuu) \right) 
+ i\, \Tr \left([E^{0}\cdot\Gamma,\czu]\right) \Big]. 
\end{align}

Finally, for greater values of $\kappa$, 
using 
(\ref{polarFiel11+}),
(\ref{polarFiel}),
(\ref{polarCoh12+}),
(\ref{polarCoh}),
(\ref{polarPop}), 
we have
\begin{align} \label{polarFiel22+}
\forall\kappa>1, \quad \fiel^{2,\kappa}=0.
\end{align}

\paragraph{Equations on the coherences.} 

When $\kappa\neq1$, we may invert directly
\begin{equation} \label{polarCoh2neq1}
\coh^{2,\kappa} = -(i\,\Og-\gamma\kappa-\vf\cdot\dthu)^{-1} 
\Big( 
\dT \coh^{1,\kappa} + \dt \coh^{0,\kappa} 
- i\, [E^0\cdot\Gamma,\coh^1+\pop^1]_{\rm od}^\kappa 
- i\, [E^1\cdot\Gamma,\coh^0+\pop^0]_{\rm od}^\kappa 
\Big) ,
\end{equation}
which vanishes for $\kappa>2$, thanks to
(\ref{polarCoh12+}),
(\ref{polarCoh}),
(\ref{polarPop1}),
(\ref{polarPop}).

For $\kappa=1$, using (\ref{polarCoh}), (\ref{polarPop}),
(\ref{polarCoh12+}), (\ref{polarPop1}) and (\ref{polarFiel11+}), 
nonresonant triples $(m,n,\alpha)$ lead to the similar formula,
\begin{align}
\label{polarCoh21}
\forall(m,n,\alpha)\notin\res, 
\quad
\coh^{2,1}_{m,n,\alpha} = 
& -(i\,\omega(m,n)+\gamma-ik\cdot\alpha_1)^{-1} 
\Big( 
 \dT \coh^{1,1}_{m,n,\alpha} 
- i\, [E^{0}\cdot\Gamma,\cuu+\puu]_{m,n,\alpha}
\\
\nonumber
& \quad
- i\, [E^{1,0}\cdot\Gamma,\czu]_{m,n,\alpha}
- i\, [E^{1,1}\cdot\Gamma,{\bf N}^0]_{m,n,\alpha} \Big) ,
\end{align}
whereas for resonant triples we obtain
\begin{align} \label{intermtranspCoh11}
\nonumber
\forall(m,n,\alpha)\in\res, 
\qquad
\dT \coh^{1,1}_{m,n,\alpha} =  
& - \dt \coh^{0,1}_{m,n,\alpha}  
 + i\, [E^{0}\cdot\Gamma,\cuu+\puu]_{m,n,\alpha} 
 \\
&
+ i\, [E^{1,0}\cdot\Gamma,\czu]_{m,n,\alpha} 
+ i\, [E^{1,1}\cdot\Gamma,{\bf N}^0]_{m,n,\alpha}.
\end{align}

\paragraph{Equations on the population.} 

For $\kappa=0$, using
(\ref{polarFiel}), 
(\ref{polarCoh}), 
(\ref{polarPop}), 
(\ref{polarCoh12+}) and 
(\ref{polarPop1}), we have
\begin{align} \label{intermtranspPop10}
&
\forall \alpha\in\cz\cup\{0\}, \qquad
\dT \pop^{1,0}_\alpha = -\dt \pop^{0}_\alpha + i\, [ E^{0}\cdot\Gamma,\cuz ]_{{\rm d},\alpha} 
+ W\op \pop^{0}_\alpha ,
\\
\label{polarPop20}
& \forall\alpha\notin\cz\cup\{0\}, \qquad
\pop^{2,0}_\alpha = (i\, \vf\cdot\alpha_1)^{-1}
\Big( \dT \pop^{1,0}_\alpha + \dt \pop^{0}_\alpha
- i\, [ E^{0}\cdot\Gamma,\cuz ]_{\rm d,\alpha} 
- W\op \pop^{0}_\alpha \Big).
\end{align}

Similarly, $\kappa>0$ leads to a polarization relation (using (\ref{polarPop})), namely
\begin{equation} \label{polarPop21+}
\forall\kappa>0, \quad  
\pop^{2,\kappa}= (\gamma\kappa+\vf\cdot\dthu)^{-1}
\Big( \dT \pop^{1,\kappa} 
- i\, [ E^0\cdot\Gamma,\coh^1 ]_{\rm d}^\kappa 
- i\, [ E^1\cdot\Gamma,\coh^0 ]_{\rm d}^\kappa \Big) ,
\end{equation}
which vanishes as soon as $\kappa>2$.

\subsection{Profile equations, intermediate scale analysis}
\label{SecInterm}

The next step in the analysis consists in obtaining a closed system
determining the first profiles $\fiel^0$, $\pop^0$ and $\coh^0$. To 
achieve this, and in order to ensure
consistancy of the Ansatz, we need impose that the corrector
terms $\fiel^1$, $\pop^1$ and $\coh^1$ are small compared with
the first profiles $\fiel^0$, $\pop^0$ and $\coh^0$.
Concerning the field and the populations, this means we need impose
$T$-sublinearity of the
correctors $\fuz$ and $\puz$ (as in \cite{Hun88}, 
\cite{JMR98}, \cite{Lan98}) while solving
equations~\eqref{intermtranspFielMoy10},
\eqref{intermtranspFielOsc10},  \eqref{intermtranspPop10}, and we shall prescribe the following requirement
\begin{equation} \label{sublin}
\forall \sigma \in \N^{3+2d}, \quad \frac{1}{T}
  \sup_{t\in[0,t_\star]} \| \partial^\sigma_{x,y,z,\theta} 
  (\fuz,\,\puz)\|_{L^2} \tendlorsque{T}{+\infty} 0.
\end{equation}
Concerning coherences, {\em i.e.} while solving
\eqref{intermtranspCoh11}, we cannot impose the same constraint 
on $\cuu$. However, and  as explained before, this corrector produces in the 
approximate solution a term 
$\rep\cuu\varTt_{|_{T=t/\rep}}e^{-\gamma t/\ep},$
so that the possible growth in $T$ of $\cuu$ is eventually compensated by the factor 
$\exp(-\gamma t/\ep)$. This is proved in
section~\ref{SecExp}.

\subsubsection{Analysis of fields and populations: average 
operators} \label{SecAvOp}

The key observation for solving equations~
\eqref{intermtranspFielMoy10}, \eqref{intermtranspFielOsc10},
and \eqref{intermtranspPop10} on
${\bf u}^{1,0}_0$, $\Pi {\bf u}^{1,0}_{\rm osc}$, and ${\bf N}^{1,0}_0$,  
 {\em while keeping}
$T$-sublinear solutions 
${\bf u}^{1,0}_0$, $\Pi {\bf u}^{1,0}_{\rm osc}$, and ${\bf N}^{1,0}_0$,  
is that the source terms in these equations
 have 
a precise stucture  in terms of propagation at the intermediate 
scale $T$.

\medskip

Let us make our point precise. Our analysis is in three steps.

\medskip

Firstly,
 equations
\eqref{intermtranspFielMoy10}, \eqref{intermtranspFielOsc10},
and \eqref{intermtranspPop10} are of the form
$$
\partial_T \Pi {\bf u}^{1,0}_{\rm osc}=\cdots, \qquad
\partial_T {\bf N}^{1,0}_0=\cdots, \qquad
\partial_T {\bf u}^{1,0}_0+M_2(0,\partial_y,\partial_z) {\bf u}^{1,0}_0 =\cdots,
$$
where the right-hand sides
 only depend on the lower order terms ${\bf u}^{0}$ and ${\bf N}^{0}$.
Besides, for $(\eta,\zeta)\in\R^2\nz$, the symmetric 
matrix $M_2(0,\eta,\zeta)$ has the spectral decomposition 
\begin{align}
M_2(0,\eta,\zeta)=\sum_{k=0,+,-} \lambda_k(\eta,\zeta) 
p_k(\eta,\zeta),
\end{align}
where the eigenvalues $\lambda_k$  are smooth on 
$\R^2\nz$ and homogeneous of degree 1, with values 
\begin{align}
\lambda_0(\eta,\zeta)=0,
\quad
\lambda_+(\eta,\zeta)=
\sqrt{\eta^2+\zeta^2},
\quad
\lambda_-(\eta,\zeta)=
-\sqrt{\eta^2+\zeta^2},
\end{align}
and the projectors $p_k$ are smooth on 
$\R^2\nz$ and homogeneous of degree 0, with the following values 
\begin{align*}
\quad
&
p_0
\text{ is the
orthogonal projector onto }
\mbox{\rm Span}\left( \begin{pmatrix} 0 \\ \mathcal{Z} \end{pmatrix} ,
  \begin{pmatrix} \mathcal{Z} \\ 0 \end{pmatrix}\right),
\\
  &
  p_{\pm}
\text{ is the
orthogonal projector onto }
  \mbox{\rm Span}\left( \begin{pmatrix} \pm |\mathcal{Z}|e_x \\
    \mathcal{Z}^\perp \end{pmatrix} , 
  \begin{pmatrix} \mathcal{Z}^\perp \\ \pm |\mathcal{Z}|e_x 
  \end{pmatrix}\right),
\\
  &
  \text{where } \quad
  \displaystyle \mathcal{Z} = \begin{pmatrix} 0 \\ \eta \\ \zeta
\end{pmatrix}, 
\quad \mathcal{Z}^\perp = \begin{pmatrix} 0 \\ -\zeta \\ \eta
\end{pmatrix}, \quad 
e_x = \begin{pmatrix} 1 \\ 0 \\ 0 \end{pmatrix} .
\end{align*}
Using Fourier transform, this allows to state the
\begin{lemme}
\label{LemStructInterm}
Take $u^{\rm in} \in H^\infty(\R^2,\C^6)$.

\noindent
Then, the u\-ni\-que so\-lu\-tion~$u \in \mathcal{C}(\R,H^\infty)$
to the Cauchy problem
$$
\partial_T u + M_2(0,\partial_y,\partial_z) u  = 0 , \qquad u_{|_{T=0}} = u^{\rm in},
$$
is given by 
$$
u = \sum_{k=0,+,-} u_k ,
$$
where each $u_k  = \pk u$ is characterized by
$$
(\dT + i\lambda_k(\Dy,\Dz)) u_k = 0 , \qquad
{u_k}_{|_{T=0}} = \pk u^{\rm in}.
$$ 
\end{lemme}
In passing,  Lemma \ref{LemStructInterm}
implies equation \eqref{intermtranspFielMoy0} induces for 
$\fiel^0_0$ the splitting 
$$
\fiel^0_0 = \fiel^0_{0,0} + \fiel^0_{0,+} + \fiel^0_{0,-} ,
$$
with
\begin{equation} \label{intermtranspFielMoy0'}
 \fiel^0_{0,k}=\pk\fiel^0_{0},\qquad
(\dT + i\lambda_k(\Dy,\Dz))\fiel^0_{0,k}=0 
\quad (k=0,+,-).
\end{equation}
Lemma \ref{LemStructInterm} also implies that equations
\eqref{intermtranspFielMoy10}, \eqref{intermtranspFielOsc10},
and \eqref{intermtranspPop10} have the form
\begin{align*}
&
\qquad\qquad
\left( \partial_T + i \lambda_k(D_y,D_z)\right) \, {\bf u}^{1,0}_{0,k}=\cdots \qquad
(k=0,+,-),\\
&
\left( \partial_T + i \lambda_0(D_y,D_z)\right) \, \Pi {\bf u}^{1,0}_{\rm osc}=\cdots,
\qquad
\left( \partial_T + i \lambda_0(D_y,D_z)\right) \, {\bf N}^{1,0}_0=\cdots,
\end{align*}
where the various right-hand-sides only depend on ${\bf u}^{0}$ and ${\bf N}^{0}$.

\medskip

With this observation in mind, the next step consists in analyzing
the above equations using  the average operators introduced in \cite{Lan98}.
They allow to describe nonlinear interactions between the various modes $0,+,-$ in the
equations at hand.

For each $k=0,+,-$, and
$u\in\mathcal{C}([0,t_\star]_t \times
\R_T,H^\infty(\R^3_{x,y,z}\times\T^{2d}_\theta))$, we define 
(omitting the dependence upon $t$, $x$ and $\theta$)  
$$
G_k^S u (T,y,z) := \frac{1}{S} \int_0^S \mathcal{F}^{-1} 
\left( e^{i \, s\lambda_k(\eta,\zeta)} \hat{u}(T+s,\eta,\zeta) 
\right) {\rm d}s,
$$ 
with $\mathcal{F}$ the Fourier transform in variables $y,z$. We also define the limit
(if it exists), 
$$
G_k u (T,y,z) := \lim_{S\to+\infty} G_k^S u (T,y,z).
$$
The average operator $G_k$ performs the average along the
bicharacteristic curves of the operator $(\dT + i\lambda_k(\Dy,\Dz))$. Naturally,
$G_0$ coincides with the usual 
average with respect to $T$, due to $\lambda_0=0$.

The following properties of the average operators $G_k$ are useful.
\begin{prop}[{\bf borrowed from \cite{Lan98}}]
\label{lannes}
Let $k \in \{0,+,-\}$. 

\vspace{0.2cm}

\noindent
(i) If $u\in\mathcal{C}([0,t_\star]_t \times 
\R_T,H^\infty(\R^3_{x,y,z}\times\T^{2d}_\theta))$ satisfies\footnote{Here, solutions to equations of the form 
$(\dT + i\lambda_k(\Dy,\Dz)) u = f$ are meant in the mild sense, 
see remark \ref{mildsol}.} 
$(\dT + i\lambda_k(\Dy,\Dz)) u = 0$, then $G_k u = u$. 

\vspace{0.2cm}

\noindent
(ii) If $f\in\mathcal{C}([0,t_\star]_t \times
\R_T,H^\infty(\R^3_{x,y,z}\times\T^{2d}_\theta))$ satisfies $G_k f = 0$, then any solution $u$ to $(\dT + i\lambda_k(\Dy,\Dz)) u = f$ is 
$T$-sublinear. 

\vspace{0.2cm}

\noindent
(iii) If $u\in\mathcal{C}^1([0,t_\star]_t \times
\R_T,H^\infty(\R^3_{x,y,z}\times\T^{2d}_\theta))$ is 
$T$-sublinear, 
then we have $G_k (\dT + i\lambda_k(\Dy,\Dz)) u = 0$. 

\vspace{0.2cm}

\noindent
(iv) Let a collection $(u_\ell)_{\,0\leq \ell \leq L} \subset
\mathcal{C}([0,t_\star]_t \times
\R_T,H^\infty(\R^3_{x,y,z}\times\T^{2d}_\theta))$ satisfy 
$(\dT + i\lambda_{k_\ell}(\Dy,\Dz)) u_\ell = 0$ with $k_\ell\in\{0,+,-\}$ for 
all $\ell$, and set $u := u_0 \cdots u_L$.

Then, if $k_\ell = k$ for 
all $\ell$, we have $G_k u = u$, else, if $k_\ell\neq k$ for some index $\ell$, we have $G_k u = 0$. 
\end{prop}

\medskip

Our last step consists in applying all above considerations. Indeed, Proposition \ref{lannes} (ii)-(iii) asserts
that equation $(\dT + i\lambda_{k}(\Dy,\Dz)) u = f$ possesses a $T$-sublinear solution if and only if $G_k f=0$, while Proposition \ref{lannes} (iv)  allows to explicitely compute $G_k f$ when $f$ is a product of solutions to $(\dT + i\lambda_{k_\ell}(\Dy,\Dz)) u_\ell = 0$.

Concerning equation \eqref{intermtranspFielMoy10},
with the wave structure given by \eqref{intermtranspFielMoy0}, we conclude that
equation \eqref{intermtranspFielMoy10} possesses a $T$-sublinear solution ${\bf u}^{1,0}_0$
if and only if
$$
\pk \mun \pk \fiel^0_{0,k} = 0 \quad (k=0,+,-).
$$
These three systems in fact reduce to the trivial condition
\begin{equation} \label{slowtranspFielMoy0}
\dt \fiel^0_{0,0} = 0, 
\quad
\dt \fiel^0_{0,+} = 0, 
\quad
\dt \fiel^0_{0,-} = 0, 
\end{equation}
thanks to the following Lemma whose proof is a straightforward computation. 
\begin{lemme} 
Take $k\in\{0,+,-\}$ and $(\eta,\zeta)\in\R^2\nz$. Then, we have
$$p_k(\eta,\zeta) A_x p_k(\eta,\zeta) = 0.$$
\end{lemme}

Concerning
equation
\eqref{intermtranspFielOsc10}, the right-hand-side of this equation only involves (products of) solutions
to $(\dT + i\lambda_{k_\ell}(\Dy,\Dz)) u_\ell = 0$, amongst which only the parts ${\bf u}^{0}_{\rm osc}$,
${\bf u}^{0}_{0,0}$ and ${\bf N}^{0}$ are associated with the characteristic speed
$\lambda_0=0$, thanks to relations
\eqref{intermtranspFielOsc0}
and \eqref{intermtranspPop0}, and thanks to the definition of
${\bf u}^{0}_{0,0}$ in \eqref{intermtranspFielMoy0'}.
Therefore, equation \eqref{intermtranspFielOsc10}
possesses a $T$-sublinear solution ${\bf u}^{1,0}_{\rm osc}$
if and only if
\begin{align}
\label{slowtranspFielOsc0}
&
(\dt+
\vDth\,\dx)\,\fzo - i \, \aDth \Delta_{y,z} \fzo = \\
\nonumber
&
\qquad
\Pi \Big( 0 , i\,\Tr 
  \left( \Gamma (i\, \Og (i\,\Og-\vf\cdot\dthu)^{-1} -1) 
  \left[ \left(E^{0}_{0,0}+E^{0}_{\rm osc}\right)
  \cdot\Gamma,\pop^{0}\right]_{\rm osc} \right) \Big).
\end{align}
This equation may be transformed  further. Indeed, when $\alpha\in\cz$, we know from Lemma 
\ref{LemVg} and Lemma \ref{LemDiff} that $\val=0$ and $\aal=0$. Besides, we already know that
the populations $\pop^{0}$ only carry temporal oscillations, according to \eqref{polarPop}.
As a consequence, we recover by a direct computation 
that equation \eqref{slowtranspFielOsc0} implies
\begin{align}
\label{slowtranspFielOsc0-bis}
&
\forall \alpha\in\cz, \quad
\dt {\bf u}^{0}_\alpha=0.
\qquad
(\text{{\em i.e. } }
\dt {\bf u}^0_{\rm space}=0).
\end{align}
(This comes from the fact  that the factor $i\, \Og (i\,\Og-\vf\cdot\dthu)^{-1} -1)$ vanishes when acting on
a frequency $\alpha\in\cz$). When $\alpha\in\mathcal{C}_{\pm}$ at variance, equation
\eqref{slowtranspFielOsc0} provides, using Lemma \ref{LemVg} and Lemma \ref{LemDiff}, the relation
\begin{align}
\label{slowtranspFielOsc0-ter}
\forall \alpha\in\mathcal{C}_{\pm},\qquad
(\dt\pm\dx)\,{\bf u}^0_\alpha \mp 
& i \, \frac{1}{2 k \cdot \alpha_0} \Delta_{y,z} {\bf u}^0_\alpha = \\
\nonumber
&
\Pi \Big( 0 , i\,\Tr 
  \left( \Gamma ( \Og (\,\Og-k\cdot \alpha_1)^{-1} -1) 
  \left[ \left(E^{0}_{0,0}+E^{0}_{\rm osc}\right)
  \cdot\Gamma,\pop^{0}\right]_{\alpha} \right) \Big).
\end{align}

Lastly, concerning equation
\eqref{intermtranspPop10}, the right-hand-side involves in the similar fashion
only products of the profiles
${\bf N}^{0}_\alpha$, $E^{0}$, and ${\bf C}^{1,0}$, which, considering the relation
\eqref{polarCoh10}, reduces to products of the profiles
${\bf N}^{0}$ and $E^{0}$.
Amongst these profiles, only the parts  ${\bf N}^{0}$, $E^{0}_{\rm osc}$, and $E^{0}_{0,0}$
are associated with the characteristic speed $\lambda_0=0$. 
Therefore, equation \eqref{intermtranspPop10}
possesses a $T$-sublinear solution ${\bf N}^{1,0}_\alpha$
if and only if (here we plug relation \eqref{polarCoh10})
\begin{align}
\label{slowtranspPop0}
\nonumber
& \forall\alpha\in\cz\cup\{0\},
\\
& \quad
\dt  \pop^0_\alpha = W\op \pop^0_\alpha
 -\, \left[
 \left( E^{0}_{0,0}+ E^{0}_{\rm osc}\right)
 \cdot\Gamma,
 (i\,\Og-\vf\cdot\dthu)^{-1} 
\left[
\left( 
E^{0}_{0,0}+
E^{0}_{\rm osc}\right)
\cdot\Gamma, \pop^{0}
\right]_{\rm od}
\right]_{{\rm d},\alpha}.
\end{align}

\medskip

Eventually, we have now obtained the set of equations
\eqref{slowtranspFielMoy0}, \eqref{slowtranspFielOsc0-bis},
\eqref{slowtranspFielOsc0-ter}, and
\eqref{slowtranspPop0} as a set of necessary and sufficient conditions to be able to find $T$-sublinear solutions to \eqref{intermtranspFielMoy10}, \eqref{intermtranspFielOsc10}, and
\eqref{intermtranspPop10}, respectively. This completes our effort in finding a system that 
completely determines the dominant profiles ${\bf u}^0$ and
${\bf N}^{0}$. 
Note in passing that equation \eqref{slowtranspFielOsc0-ter} in fact reduces to
\begin{align}
\label{slowtranspFielOsc0-qua}
\forall \alpha\in\mathcal{C}_{\pm},\qquad
(\dt\pm\dx)\,{\bf u}^0_\alpha \mp 
& i \, \frac{1}{2 k \cdot \alpha_0} \Delta_{y,z} {\bf u}^0_\alpha = \\
\nonumber
&
\qquad
\Pi \Big( 0 , i\,\Tr 
  \left( \Gamma ( \Og (\,\Og-k\cdot \alpha_1)^{-1} -1) 
  \left[ E^{0}_{\rm time} \cdot\Gamma,\pop^{0}\right]_{\alpha} \right) \Big),
\end{align}
due to the fact that ${\bf N}^0_\alpha=0$ whenever $\alpha\notin\cz\cup\{0\}$.

In conclusion, we have recovered equations 
\eqref{coco1}, 
\eqref{bahbouh}, 
\eqref{nltran}, 
\eqref{22}, 
\eqref{23}, 
\eqref{dtnz}, 
\eqref{rate} that had been announced in the introductory part of this paper.

\subsubsection{Analysis of coherences: exponential growth} 
\label{SecExp}

The above analysis is not possible in the case of 
\eqref{intermtranspCoh11} :  equation
\eqref{intermtranspCoh0}, ruling the evolution of $\czu$ as a function of $T$, is not a constant 
coefficient system. Instead, we consider it as a non-autonomous 
system of linear ODE's, parametrized by $t$, $x$, $y$, $z$, 
$\theta$, with smooth and bounded coefficients. This point  of view provides the simple

\begin{lemme} \label{LemExp}
Let
$(\fiel^0,\pop^0)\in\mathcal{C}([0,t_\star]_t\times\R_T,
H^\infty(\R^3_{x,y,z}\times\T^{2d}_\theta))$ be a mild solution to 
\eqref{intermtranspFielMoy0}, \eqref{intermtranspFielOsc0},  
\eqref{intermtranspPop0},
\eqref{slowtranspFielMoy0}, \eqref{slowtranspFielOsc0}, 
\eqref{slowtranspPop0}. Then, 
\begin{itemize}
\item[(i)] $(\fiel^0,\pop^0)\in\mathcal{C}^\infty_{\rm b}([0,t_\star]_t\times\R_T,
H^\infty(\R^3_{x,y,z}\times\T^{2d}_\theta))$ is smooth and bounded uniformly with respect to all variables, as well as all its 
derivatives. 
\item[(ii)] associated with these values of $(\fiel^0,\pop^0)$, take any solution $\czu$ to \eqref{intermtranspCoh0} 
with an initial data belonging to $\mathcal{C}^\infty\left( [0,t_\star]_t,H^\infty(\R^3_{x,y,z}\times\T^{2d}_\theta)\right)$. Then,
for any multi-index $\mu\in\N^{5+2d}$, there are constants 
$K_1,K_2>0$ such that, uniformly on $[0,t_\star]_t\times\R_T
\times\R^3_{x,y,z}\times\T^{2d}_\theta$, we have
$$\left|\partial_{t,T,x,y,z,\theta}^\mu\czu(T)\right|\leq 
K_1e^{K_2T}.$$
\end{itemize}
\end{lemme}

\begin{rem}
\label{mildsol}
All equations that are referred to here are either of the form $\dt v=f(v,Dv,D^2v)$, or of the form
$\dT v=f(v,Dv,D^2v)$, where $v$ is assumed to have $H^\infty$ smoothness in $(x,y,z,\theta)$, where the symbol $D$ means differentiation with respect
to $(x,y,z,\theta)$, and $f$ is a possibly non-linear function that depends on the equation. The notion of mild solution is then the standard one :
we mean a solution to the integral equation $v(t)=v|_{t=0}+\int_0^t f(v,Dv,D^2v)(s) \, {\rm d} s$ or  to
$v(T)=v|_{T=0}+\int_0^T f(v,Dv,D^2v)(S) \, {\rm d} S$.
\end{rem}

\noindent
{\bf Proof.} Boundedness and smoothness of $(\fiel^0,\pop^0)$ with respect to $(x,y,z,\theta)$ is obvious.
Concerning the variables $t$ and $T$, we exploit the fact that $t$ belongs to a compact set.
We also exploit  the structure
of the relevant equations. More precisely, $\fiel^0_0$ satisfies $\left(\dT+M_2(0,\dy,\dz)\right) \fiel_0^0=0$ and $\dt \fiel^0_0=0$. This, together with Lemma
\ref{LemStructInterm}, provides boundedness and smoothness of $\fiel^0_0$ with respect to $t$ and $T$.
Next, $\pop^0$ satisfies $\dT \pop^0=0$ together with equation \eqref{slowtranspPop0}, an equation
of the form $\dt\pop^0=f(\fiel^0,\pop^0)$ where
$f$ is smooth.
This provides boundedness and smoothness of $\pop^0$ with respect to $t$ and $T$. We have used here the 
Diophantine Hypothesis  \ref{HypDioph} to conclude that $f$ is smooth.
Lastly, $\fiel^0_{\rm osc}$ satisfies $\dT \fiel^0_{\rm osc}=0$ together with equation
\eqref{slowtranspFielOsc0}, an equation of the form
$\dt \fiel^0_{\rm osc}=f\left(\fiel^0_{\rm osc},D_\theta \fiel^0_{\rm osc},D^2_\theta\fiel^0_{\rm osc},\fiel^0_0,\pop^0 \right)$,
where $f$ is smooth.
This provides boundedness and smoothness of $\fiel^0_{\rm osc}$ with respect to $t$ and $T$.

Point \emph{(ii)} now comes from the  Gronwall lemma, together with the fact that equation \eqref{intermtranspCoh0}  has the form
$\dT {\bf C}^{0,1}=f\left(E^0,{\bf C}^{0,1}\right)$ where $f$ is bilinear. 
\fin

\begin{rem}
Under the assumptions of the above Lemma, pushing the analysis further, we may consider a solution ${\bf C}^{1,1}$ to
 \eqref{intermtranspCoh11}, and look for the available estimates on ${\bf C}^{1,1}$. This function 
 is associated with a solution ${\bf N}^{1,1}$ to \eqref{polarPop1}. In view of the above Lemma, the latter obviously satisfies an exponential bound
 of the form 
 $\left|\partial_{t,T,x,y,z,\theta}^\mu {\bf N}^{1,1}(T)\right|\leq 
K_1 \, e^{K_2T},$ where $K_1$ and $K_2$ depend on $\sigma$ but not on $T$.
Hence ${\bf C}^{1,1}$ satisfies similarly
$\left|\partial_{t,T,x,y,z,\theta}^\mu {\bf C}^{1,1}(T)\right|\leq 
K_1 \, e^{K_2T}.$ As a consequence, we recover 
\begin{equation} \label{estimCoh11}
\left| \rep \cuu \left(t,x,y,z,t/\rep,\theta \right) 
e^{-\gamma t/\ep} \right| \leq  K_1 \rep 
e^{(K_2-\gamma/\rep)\,t/\rep} = \mathcal{O}(\rep).
\end{equation}
\end{rem}

\subsection{Solving the profile equations} \label{SecExist}

\subsubsection{Computing the dominant profile ${\bf U}^0$.}

The first profile $\sol^0=(\fiel^0,\coh^0,\pop^0)$ is constrained by the polarization 
conditions 
\eqref{polarFiel}, \eqref{polarCoh}, \eqref{polarPop}. Besides,  
its components are propagated in various ways. The average 
$\fiel^0_0$  does not depend on the slow time 
$t$ (equation \eqref{slowtranspFielMoy0}), and it satisfies the linear hyperbolic equation
$\left(\dT + M_2(0,\dy,\dz)\right) \fiel^0_0=0$ with respect to
the intermediate time $T$ (equation \eqref{intermtranspFielMoy0'}). The oscillating part $\fzo$, and the populations $\pop^0$, do not depend on $T$ 
(equations \eqref{intermtranspFielOsc0}, \eqref{intermtranspPop0}), 
and they satisfy   nonlinear evolution equations with respect to time 
$t$ (equations \eqref{slowtranspFielOsc0}, \eqref{slowtranspPop0}). Coherences $\coh^0$ are only constrained
to satisfy  the nonlinear 
ODE \eqref{intermtranspCoh0} in time $T$, and the slow time  $t$ only acts as a parameter here.

We are now in position to state the

\begin{theo} \label{ThExistPro} Let $s>(3+2d)/2$. 
Take a function\footnote{The reader should be cautious about the fact that the initial datum only depends on $\theta_0$, and not on $\theta_1$. Recall that eventually $\theta_0$ will be replaced by $kx/\ep$ while $\theta_1$ takes the value $-kt/\ep$.}
$$
\underline{\bf U}^0 \equiv \underline{\bf U}^0(x,y,z,\theta_0),
$$
which belongs to
$H^s(\R^3_{x,y,z}\times \T^d_{\theta_0})$.
Assume that the component $\underline{\coh}^0$ of $\underline{\bf U}^0$
satisfies the following polarization condition
\begin{equation}
\label{polarCohbibis}
\forall (m,n) \in \Z^{2d} \text{ satisfying } (m,n,\alpha) \notin \mathcal{R}(k)
\text{ for all }Ê\alpha\in \Z^{2d}, 
\quad
\text{ we have } \quad \underline{\coh}^0_{m,n}=0.
\end{equation}
Then, 
there is $t_\star>0$ and a unique function 
$$
{\bf U}^0 = ({\bf u}^0,\czu e^{-\sigma},{\bf N}^0)
\in \mathcal{C}([0,t_\star]_t\times[0,+\infty[_T\times[0,+\infty[_\sigma,
H^s(\R^3_{x,y,z}\times\T^{2d}_{(\theta_0,\theta_1)})),
$$
solution to
the polarization conditions \eqref{polarFiel}, \eqref{polarCoh}, \eqref{polarPop}, 
and which satisfies
equations \eqref{intermtranspFielMoy0'}, \eqref{slowtranspFielMoy0}
for the average field ${\bf u}^0_0$, equations 
\eqref{intermtranspFielOsc0}, \eqref{slowtranspFielOsc0}
for the oscillatory field ${\bf u}^0_{\rm osc}$, equations 
\eqref{intermtranspPop0},
 \eqref{slowtranspPop0}
for the populations ${\bf N}^0$, and equation
\eqref{intermtranspCoh0} for the coherences ${\bf C}^{0,1}$.
Uniqueness and existence is guaranteed 
provided we impose ${\bf U}^0$ 
satisfies besides the two initial constraints\footnote{Due to 
the fact that the coherences only need to satisfy an evolution equation 
in the intermediate time $T$, note that the slow time $t$
needs to be treated separately here, including in terms of initial data. Note also that the initial data is {\em not}
decomposed into modes, contrary to the solution itself.}
\begin{align}
\label{initconcon}
{\bf U}^0_{|_{t=T=\sigma=0, \theta_1=0}}=\underline{\bf U}^0,
\qquad
\text{ and } \, \czu_{|_{T=0}} \text{ is independent of } t.
\end{align}

If in addition we take an integer $\ell\in\N$ and assume that $s>(3+2d)/2+\ell(a+2)$, 
where the positive real $a$ is the one entering  the Diophantine Hypothesis~\ref{HypDioph},
then we recover the higher regularity
$$
{\bf U}^0 \in \mathcal{C}^\ell\left(
[0,t_\star]_t\times[0,+\infty[_T\times[0,+\infty[_\sigma,
H^{s-\ell(a+2)}\left(
\R^3_{x,y,z}\times\T^{2d}_{(\theta_0,\theta_1)}
\right)
\right). 
$$
\end{theo}

\noindent
{\bf Proof.}

{\em First step : identifying the initial values of the profiles.}

The first and most important step consists in understanding the polarization 
conditions, and how they generate the relevant initial data for the different modes. 
The data $\underline{\bf U}^0(x,y,z,\theta_0)$ is the
value of  the solution ${\bf U}^{0,\kappa}$ at $t=T=\sigma=0$ and $\theta_1=0$. 
We prove that  it defines the value of ${\bf U}^{0,\kappa}$ at $t=T=0$ for all values 
of $\theta_1\in\T^d$. 

Firstly, concerning the fields,
the polarization constraint  \eqref{polarFiel} provides, for the 
(to be defined) solution $\fiel^0$, the relations
$\fiel^0=\fzz$ and  $\Pi \fiel^0=\fiel^0$.
On the other hand, we may write down the Fourier transform of 
$\underline{\bf u}^0$ in the sole variable
$\theta_0$, together with the Fourier transform of $\fiel^0$ 
in $(\theta_0,\theta_1)$, and obtain 
\begin{equation*}
\underline{\fiel}^0(x,y,z,\theta_0) = \sum_{\beta\in\Z^d} 
\underline{\fiel}^0_\beta(x,y,z) e^{i\beta\cdot\theta_0}, \qquad
\fiel^0(t,T,\sigma,x,y,z,\theta_0,\theta_1) = \sum_{\alpha\in\mathcal{C}\cup\{0\}} 
\fiel^0_\alpha(t,T,x,y,z) e^{i\alpha\cdot(\theta_0,\theta_1)}.
\end{equation*}
Taking $t=T=\sigma=0$, $\theta_1=0$ in the second equation, 
and equating Fourier coefficients, then provides
\begin{equation*}
\forall \beta\in\Z^{d}\setminus\{0\}, \quad
\underline{\fiel}^0_\beta
=
\left(
\fiel^0_{(\beta,\beta)}
+
\fiel^0_{(\beta,-\beta)}
+
\fiel^0_{(\beta,0)}
\right)\Big|_{t=T=0},
\qquad\qquad
\underline{\fiel}^0_0
=
\fiel^0_{(0,0)}\Big|_{t=T=0},
\end{equation*}
where 
$\fiel^0_{(\beta,\beta)}=\pi^+ \, \fiel^0_{(\beta,\beta)}$, and
$\fiel^0_{(\beta,-\beta)}=\pi^- \, \fiel^0_{(\beta,-\beta)}$, and
$\fiel^0_{(\beta,0)}=\pi^0 \, \fiel^0_{(\beta,0)}$,
thanks to the polarization conditions and using the value of $\Pi$ (see \eqref{pipi}).
Hence, applying successively the three operators $\pi^+$, $\pi^-$, $\pi^0$ on both sides
of the first equality, and using the obvious orthogonality
relations $\pi^+\pi^-=0$ and so on, we recover
the necessary relations
\begin{align*}
&
\fiel^0_{(0,0)}\Big|_{t=T=0}=
\underline{\fiel}^0_0, \qquad \text{ and, when } \beta\neq 0,
\\
&\qquad
\fiel^0_{(\beta,\beta)}\Big|_{t=T=0}
=
\pi^+ \, \underline{\fiel}^0_\beta,
\quad
\fiel^0_{(\beta,-\beta)}\Big|_{t=T=0}
=
\pi^- \, \underline{\fiel}^0_\beta,
\quad
\fiel^0_{(\beta,0)}\Big|_{t=T=0}
=
\pi^0 \, \underline{\fiel}^0_\beta.
\end{align*}
This terminates the analysis of the initial conditions for the fields.

Secondly, concerning the populations, the polarization constraint 
\eqref{polarPop} provides the simpler relations 
${\bf N}^0={\bf N}^{0,0}$ and ${\bf N}^0_\alpha=0$ whenever 
$\alpha\notin\cz\cup\{0\}$. Writing down the Fourier transforms of 
$\underline{\bf N}^0$ and ${\bf N}^0$ as in the previous paragraph 
then provides
$
\displaystyle
\underline{\bf N}^0(x,y,z,\theta_0) = \sum_{\beta\in\Z^d} 
\underline{\bf N}^0_\beta(x,y,z) e^{i\beta\cdot\theta_0},
$
and
$
{\bf N}^0(t,T,\sigma,x,y,z,\theta_0, \theta_1) = \sum_{\beta\in\Z^d} 
{\bf N}^0_\beta(t,T,x,y,z) e^{i\beta\cdot\theta_0}.
$
Specifying  $t=T=\sigma=0$, $\theta_1=0$, and equating Fourier 
coefficients, gives the simple value
\begin{align*}
&
\forall \beta\in\Z^d, \quad
{\bf N}^0_\beta\Big|_{t=T=0}=
\underline{\bf N}^0_\beta.
\end{align*}
This terminates the analysis of the initial conditions for 
the populations.

Lastly, concerning the coherences, the polarization conditions 
\eqref{polarCoh} asserts ${\bf C}^{0,1}_{m,n,\alpha}=0$ whenever
$(m,n,\alpha)\notin \mathcal{R}(k)$. Writing down the natural 
expansions of ${\bf C}^0={\bf C}^{0,1} \, e^{-\sigma}$ and
$\underline{\bf C}^0$, we recover for each value of $m$ and $n$ 
the identities
$
\underline{\bf C}^0_{m,n}(x,y,z,\theta_0) = \sum_{\beta\in\Z^d} 
\underline{\bf C}^0_{m,n,\beta}(x,y,z) e^{i\beta\cdot\theta_0}
$,
together with
$
{\bf C}^0_{m,n}(t,T,\sigma,x,y,z,\theta_0,\theta_1) =\left( \sum_{\alpha\in\Z^{2d}} 
{\bf C}^{0,1}_{m,n,\alpha}(t,T,x,y,z) e^{i\alpha\cdot(\theta_0,\theta_1)}\right) \, e^{-\sigma}.
$
Specifying $t=T=\sigma=0$, $\theta_1=0$, 
gives using the polarization conditions
\begin{align}
\label{2be}
&
\forall \beta\in\Z^d, \quad
\underline{\bf C}^0_{m,n,\beta}
=
\sum_{\alpha_1\in\Z^d}
{\bf C}^{0,1}_{m,n,(\beta,\alpha_1)}.
\end{align}
Now, two cases occur, depending on the value of $(m,n)$.

\indent
$\bullet$ First case : for any $\alpha\in\Z^{2d}$, we have $(m,n,\alpha)\notin \mathcal{R}(k)$.

In that case, equation \eqref{2be} is automatically satisfied since the right-hand-side vanishes
due to the polarization constraint on ${\bf C}^{0,1}$, while the left-hand-side
is assumed to vanish thanks to the additional constraint
\eqref{polarCohbibis} we have set on the initial function $\underline{\bf C}^0$.

\indent
$\bullet$ Second case : there exists an $\alpha\in\Z^{2d}$, 
such that $(m,n,\alpha)\in \mathcal{R}(k)$.
 
In that case, the solution $\widetilde\alpha_1\in\Z^d$ to the 
equation $\omega(m,n)=k \cdot \widetilde\alpha_1$ is unique, 
thanks to the $\Q$-independence of the coordinates of the wave-vector 
$k$. 
Recall that this equation defines the set $ \mathcal{R}(k)$.
Therefore, given any $\beta\in \Z^d$, equation \eqref{2be} reduces to
\begin{align*}
{\bf C}^{0,1}_{m,n,(\beta,\widetilde\alpha_1)}
=
\underline{\bf C}^0_{m,n,\beta},
\qquad
\text{ and } \forall\alpha_1\neq\widetilde\alpha_1, \quad
{\bf C}^{0,1}_{m,n,(\beta,\alpha_1)}=0.
\end{align*}
This terminates the analysis of the initial conditions for the coherences.

\bigskip

{\em Second step : solving the evolution equations.}

We now consider the evolution problem with respect to the times $t$ and $T$, 
with the above derived initial data.

The field average ${\bf u}^0_0$ decouples and may be determined first starting from the
initial value  ${\fiel_0^0}_{|_{t=T=0}}$, by using
$(\dT+M_2(0,\dy,\dz)) {\bf u}^0_0 =0$ and $\dt {\bf u}^0_0=0$.

Next, we may solve the equations
on $({\bf u}^0_{\rm osc},{\bf N}^0)$.
The evolution with respect to $T$ is trivial since
$\dT ({\bf u}^0_{\rm osc},{\bf N}^0)=0$.
There remains to solve the coupled nonlinear system \eqref{slowtranspFielOsc0},
\eqref{slowtranspPop0}, an evolution equation in $t$. 
To do so, we use a standard iterative scheme 
(see for example \cite{AG91}), and introduce the iteration
\begin{align*} 
(\dt+\vDth\,\dx)\,{\bf u}^{0,(n+1)}_{\rm osc} 
- i & \aDth \Delta_{y,z} {\bf u}^{0,(n+1)}_{\rm osc} =
\\
& \Pi \Big( 0 , i\Tr \left( \Gamma ( \Og \,(\Og-\vf\cdot\dthz)^{-1} -1) 
[E^{0,(n)}_{\rm osc}\cdot\Gamma,\pop^{0,(n)}]_{\rm osc} \right) \Big),
\\ 
\dt  \pop^{0,(n+1)} = W\op  \pop^{0,(n+1)} 
& - \Big[ ( E^{0,(n)}_{\rm osc}+ E^0_{0,0})\cdot\Gamma ,
(i\,\Og-\vf\cdot\dthz)^{-1} \crod{\left(E^{0,(n)}_{\rm osc}+
  E^0_{0,0}\right)\cdot\Gamma}{\pop^{0,(n)}} \Big]_{\rm d,space}
\\ 
& - \Big[ ( E^{0,(n)}_{\rm osc}+ E^0_{0,0})\cdot\Gamma ,
(i\,\Og-\vf\cdot\dthz)^{-1} \crod{\left(E^{0,(n)}_{\rm osc}+
  E^0_{0,0}\right)\cdot\Gamma}{\pop^{0,(n)}} \Big]_{\rm d,0},
\end{align*}
with initial data $\left(\left(\fiel^0_{\rm osc}\right)_{|_{t=T=0}},
{\pop^0}_{|_{t=T=0}}\right)$, and an initial "guess" set to, say, 
$(\fiel^{0,(0)}_{\rm osc},\pop^{0,(0)}) = {\rm const}
= \left(\left(\fiel^0_{\rm osc}\right)_{|_{t=T=0}}, 
{\pop^0}_{|_{t=T=0}}\right)$ for any $t$ (and $T$). 
For this linearized scheme, usual energy estimates are
available in any Sobolev space $H^s$,  
to which the skew-symmetric operators
$\vDth\,\dx$ and $i \, \aDth \Delta_{y,z}$ do not contribute. 
They are
\begin{equation*}
\left\|U^{0,(m+1)}(t)-U^{0,(n+1)}(t)\right\|_{H^s} 
\leq 
\int_0^t e^{C(t-s)}
\left\|F(U^{0,(m)})(s)-F(U^{0,(n)})(s)
\right\|_{H^s} {\rm d} s,
\end{equation*}
where $F({U^0}^{(n)})$ stands for the right-hand-side in the above 
iteration.
Now, standard nonlinear tools assert that the function $F$ acts
in a locally Lipschitz fashion on $H^s=H^s(\R^3\times\T^{2d})$ 
provided $s>(3+2d)/2$.
The existence and uniqueness of a solution 
$({\bf u}^0_{\rm osc},{\bf N}^0)$ to \eqref{slowtranspFielOsc0},
\eqref{slowtranspPop0} on some time interval $[0,t_*]$ then follows.

Concerning the  coherences $\czu$, the statement of our Theorem 
imposes to choose them independent of $t$ at time $T=0$. 
This constraint is only a (pratical) way to fix the value of  
${\bf C}^{0,1}$ on the set $\{T=0\}$ (recall that ${\bf C}^{0,1}$ 
satisfies an evolution equation in time $T$ only). 
Equation \eqref{intermtranspCoh0} then asserts 
$$
\dT {\bf C}^{0,1}=i\, [E^0 \cdot \Gamma, {\bf C}^{0,1}].
$$
This equation is enough to uniquely determine ${\bf C}^{0,1}$ 
for any $t$ and $T$, starting from its known values on $\{T=0\}$.

The polarization conditions \eqref{polarFiel}, \eqref{polarCoh},
\eqref{polarPop} commute with equations 
\eqref{intermtranspFielMoy0'}, \eqref{slowtranspFielMoy0}
for the average field ${\bf u}^0_0$, with equations 
\eqref{intermtranspFielOsc0}, \eqref{slowtranspFielOsc0}
for the oscillatory field ${\bf u}^0_{\rm osc}$, with equations 
\eqref{intermtranspPop0},
 \eqref{slowtranspPop0}
for the populations ${\bf N}^0$, and with equation
\eqref{intermtranspCoh0} for the coherences ${\bf C}^{0,1}$.
Hence
by
uniqueness of the solutions to these propagation equations,
polarizations are preserved along the evolution.

There remains to study the higher regularity of the solutions we have exhibited.
The equations on the mean field, the equations on the populations,
and the equations on the coherences clearly do not induce any loss of smoothness,
{\em i.e.} provided ${\bf u}^0_{\rm osc}$ is $H^\infty$,
the derivatives
$\dt({\bf u}^0_0,{\bf C}^{0,1},{\bf N}^{0})$ and 
$\dT({\bf u}^0_0,{\bf C}^{0,1},{\bf N}^{0})$ have the same $H^s$ smoothness
as $({\bf u}^0_0,{\bf C}^{0,1},{\bf N}^{0})$. 
The loss of smoothness comes from the Schr\"odinger-like equation
\eqref{slowtranspFielOsc0} on ${\bf u}^0_{\rm osc}$. When differentiating
this equation with respect to $t$ indeed, an additional factor
$\displaystyle \aal \Delta_{y,z}= \pm \frac{1}{2\vo\cdot\alpha_0} \Delta_{y,z}$
comes into play which causes,
from the Diophantine Hypothesis \ref{HypDioph}, a loss of $s_0=a+2$ derivatives
(two derivatives are lost due to the Laplacian, while $a$ derivatives are lost due to 
small divisors, seen in Fourier
space). This argument explains why
${\bf U}^0 \in \mathcal{C}^1(H^{s-s_0})$ and, more generally, why  
${\bf U}^0 \in \mathcal{C}^\ell(H^{s-\ell \, s_0})$.
\fin

\subsubsection{Computing the correctors ${\bf U}^1$ and ${\bf U}^2$.}

The first corrector $\sol^1=({\bf u}^1,{\bf C}^1,{\bf N}^1)$ is built up in the following way.
On the one hand, we have
\begin{equation}
\label{ploum}
{\bf u}^1={\bf u}^{1,0}+{\bf u}^{1,1} \, e^{-\sigma}, \qquad
{\bf C}^1={\bf C}^{1,0}+{\bf C}^{1,1} \, e^{-\sigma}, \qquad 
{\bf N}^1=\left({\bf N}^{1,0}_0+{\bf N}^{1,0}_{\rm space}\right) 
+ {\bf N}^{1,1} \, e^{-\sigma},
\end{equation}
thanks to \eqref{polarFiel11+} , \eqref{polarCoh12+} ,
\eqref{polarPop1} .
Besides,
the contributions
$(1-\Pi) {\bf u}^{1,0}_{\rm osc}$,
${\bf u}^{1,1}$,
${\bf C}^{1,0}$,
${\bf C}^{1,1}_{m,n,\alpha} \, {\bf 1}_{(m,n,\alpha)\notin \mathcal{R}(k)}$,
${\bf N}^{1,1}$,
are prescribed as functions of the dominant profile ${\bf U}^0$ 
through the polarization conditions
\eqref{polarFiel10},
\eqref{polarFiel11}, 
\eqref{polarCoh10}, 
\eqref{polarCoh11},
\eqref{polarPop1},
respectively. They are thus known and smooth functions of ${\bf U}^0$.
Let us write these relations symbolically as 
\begin{align}
\label{bouboubou}
\left((1-\Pi) {\bf u}^{1,0}_{\rm osc},
{\bf u}^{1,1}, 
{\bf C}^{1,0},
{\bf C}^{1,1}_{m,n,\alpha} \, {\bf 1}_{(m,n,\alpha)\notin \mathcal{R}(k)},
{\bf N}^{1,1}\right)=\phi\left({\bf U}^0\right).
\end{align}
The remaining terms 
\begin{align*}
{\bf u}^{1,0}_0, \quad
\Pi {\bf u}^{1,0}_{\rm osc}, \quad 
{\bf C}^{1,1}_{m,n,\alpha} \, {\bf 1}_{(m,n,\alpha)\in \mathcal{R}(k)}, \quad
{\bf N}^{1,0}_0+{\bf N}^{1,0}_{\rm space},
\end{align*}
are then determined through linear evolution equations in $T$, namely through
\eqref{intermtranspFielMoy10},
\eqref{intermtranspFielOsc10},
\eqref{intermtranspCoh11},
\eqref{intermtranspPop10}, respectively. 
These equations read, after some simplifications,
\begin{align}
&
\nonumber
\left(\dT +M_2(0,\dy,\dz)\right) \fiel^{1,0}_0 = -A_x\dx\fiel^{0,0}_0 ,
\qquad
\dT \Pi \fuzo = 0 ,
\\
&
\nonumber
\forall\alpha\in\cz\cup\{0\}, \quad
\dT  {\bf N}^{1,0}_\alpha =
- \Big[ E^0\cdot\Gamma , (i\,\Og-\vf\cdot\dthz)^{-1} 
\crod{E^0\cdot\Gamma}{\pop^0} \Big]_{{\rm d},\alpha}
\\
& 
\nonumber
\qquad\qquad\qquad\qquad\qquad\qquad
+ \Big[ (E^0_{\rm osc} + E^0_{0,0})\cdot\Gamma ,
(i\,\Og-\vf\cdot\dthz)^{-1} \crod{(E^0_{\rm osc} +
  E^0_{0,0})\cdot\Gamma}{\pop^0} \Big]_{{\rm d},\alpha},
\\
&
\nonumber
\forall(m,n,\alpha)\in\res, \quad
\dT \coh^{1,1}_{m,n,\alpha} =  - \dt \coh^{0,1}_{m,n,\alpha}  
+ i\, [E^{0,0}\cdot\Gamma,\cuu+\puu]_{m,n,\alpha} 
\\
&
\label{intermtranspPop10'}
\qquad\qquad\qquad\qquad\qquad\qquad\qquad\qquad
+ i\, [E^{1,0}\cdot\Gamma,\czu]_{m,n,\alpha}
+ i\, [E^{1,1}\cdot\Gamma,\pzz]_{m,n,\alpha}.
\end{align}
Note that the right-hand-side of the above system only involves
known functions of ${\bf U}^0$ at this stage. In order to solve 
this system, there remains to impose as we did for $\czu$, 
\begin{align}
\label{bibibi}
\left(\fiel^{1,0}_0, \Pi\fuzo, {\bf N}^{1,0}_\alpha \, {\bf 1}_{\alpha\in\cz\cup\{0\}},
\coh^{1,1}_{m,n,\alpha} \, {\bf 1}_{(m,n,\alpha)\in \mathcal{R}(k)}\right)\Big|_{T=0}
\text{ is independent of } t.
\end{align}
This completely determines ${\bf U}^1$ as a function of ${\bf U}^0$.

The second corrector $\sol^2=({\bf u}^2,{\bf C}^2,{\bf N}^2)$ 
is built up in the following way. On the one hand, we have
\begin{equation}
\label{ploumploum}
{\bf u}^2={\bf u}^{2,0}+{\bf u}^{2,1} \, e^{-\sigma}, \quad
{\bf C}^2={\bf C}^{2,0}+{\bf C}^{2,1} \, e^{-\sigma}+{\bf C}^{2,2} 
\, e^{-2\sigma}, \quad
{\bf N}^2={\bf N}^{2,0} +{\bf N}^{2,1} \, e^{-\sigma}+{\bf N}^{2,2} 
\, e^{-2\sigma},
\end{equation}
thanks to \eqref{polarFiel22+} , \eqref{polarCoh2neq1} ,
\eqref{polarPop21+}. Besides, the contributions
$(1-\Pi) {\bf u}^{2,0}_{\rm osc}$,
${\bf u}^{2,1}$,
${\bf C}^{2,0}$,
${\bf C}^{2,1}_{m,n,\alpha} \, {\bf 1}_{(m,n,\alpha)\notin \mathcal{R}(k)}$,
${\bf C}^{2,2}$,
${\bf N}^{2,0}_{\alpha} \, {\bf 1}_{\alpha\notin\cz\cup\{0\}}$,
${\bf N}^{2,1}$,
${\bf N}^{2,2}$,
are prescribed as functions on the two first profiles ${\bf U}^0$ and ${\bf U}^1$ through the
polarization conditions
\eqref{polarFiel20},
\eqref{polarFiel21}, 
\eqref{polarCoh2neq1}, 
\eqref{polarCoh21},
\eqref{polarCoh2neq1}, 
\eqref{polarPop20},
\eqref{polarPop21+},
\eqref{polarPop21+}, 
respectively. They are thus known, smooth functions of ${\bf U}^0$ and ${\bf U}^1$.
Some parts of the corrector are free and may be chosen arbitrarily, namely
$\fiel^{2,0}_0$, 
$\Pi\fiel^{2,0}_{\rm osc}$,
$\coh^{2,1}_{m,n,\alpha} \, {\bf 1}_{(m,n,\alpha)\in \mathcal{R}(k)}$,
$\pop^{2,0}_\alpha \, {\bf 1}_{\alpha\in\cz\cup\{0\}}$. We make the most simple choice to set these contributions to zero.
Let us write all these relations symbolically as 
\begin{align}
\label{boubouboubou}
{\bf U}^2=\psi\left({\bf U}^0,{\bf U}^1\right).
\end{align}
\medskip

All these considerations, in conjunction with Theorem~\ref{ThExistPro},
lead to the
\begin{prop} \label{PropExistCorr}
For all $\ell\in\N$, there exists $\tilde{s}(\ell,a)$ such that 
the following holds.

Take an $s>\tilde{s}+(3+2d)/2$. Take an initial data
$\underline{{\bf U}}^0$ in $H^s(\R^3_{x,y,z}\times\T^d_{\theta_0})$ 
which satisfies the constraint \eqref{polarCohbibis}. 
Take the profile $\sol^0$ associated with these data 
through Theorem~\ref{ThExistPro}. 
Take two initial data 
$\underline{\sol}^1$ in $H^s(\R^3_{x,y,z}\times\T^d_{\theta_0})$ and 
$\underline{\sol}^2$ in $H^s(\R^3_{x,y,z}\times\T^d_{\theta_0})$. 
Then, there are unique correctors
$$\sol^1=
\begin{pmatrix}
\fuz+\fiel^{1,1} e^{-\sigma} \\ 
\cuz+\coh^{1,1} \, e^{-\sigma} \\
\left({\bf N}^{1,0}_0+{\bf N}^{1,0}_{\rm space}\right)+\pop^{1,1} e^{-\sigma}
\end{pmatrix}, 
\quad
\sol^2=
\begin{pmatrix}
\fiel^{2,0}+\fiel^{2,1} e^{-\sigma} \\ 
\coh^{2,0}+\coh^{2,1}  e^{-\sigma}+
\coh^{2,2} e^{-2\sigma} \\ 
\pop^{2,0}+\pop^{2,1} e^{-\sigma}+\pop^{2,2} e^{-2\sigma}
\end{pmatrix}
$$ 
in $\mathcal{C}^\ell([0,t_\star]_t\times[0,+\infty[_T\times
[0,+\infty[_\sigma,H^{s-\tilde{s}}(\R^3_{x,y,z}\times
\T^{2d}_{(\theta_0,\theta_1)}))$, 
which satisfy the constraints \eqref{bouboubou}, 
\eqref{bibibi} for the first corrector, the constraints 
\eqref{boubouboubou} for the second corrector, 
and the evolution equation \eqref{intermtranspPop10} 
for the first corrector, together with the initial constraint
$$
\sol^1\Big|_{T=t=\sigma=0,  \theta_1=0}
=
\underline{\sol}^1
\qquad
\sol^2\Big|_{T=t=\sigma=0, \theta_1=0}
=
\underline{\sol}^2.
$$
\end{prop}

\noindent
{\bf Proof.}
The initial constraint
$\sol^1\Big|_{T=t=\sigma=0, \theta_1=0}=\underline{\sol}^1$ 
reads
\begin{align*}
\left(
\fuz+\fiel^{1,1},
\cuz+\coh^{1,1},
\left({\bf N}^{1,0}_0+{\bf N}^{1,0}_{\rm space}\right)+\pop^{1,1} 
\right)\Big|_{T=t=0, \theta_1=0}
=
\underline{\sol}^1.
\end{align*}
This, together with
\eqref{bouboubou} written at time $T=t=0$,
prescribes the value of 
\begin{align*}
&
\left((1-\Pi) {\bf u}^{1,0}_{\rm osc},
{\bf u}^{1,1}, 
{\bf C}^{1,0},
{\bf C}^{1,1}_{n,m,\alpha} \, {\bf 1}_{(m,n,\alpha)\notin \mathcal{R}(k)},
{\bf N}^{1,1}\right)\Big|_{T=t=0}=\phi\left(\underline{\bf U}^0\right).
\end{align*}
Hence, taking the difference, we recover the value of 
\begin{align*}
&
\left(
{\bf u}^{1,0}_0
+
\Pi {\bf u}^{1,0}_{\rm osc},
\sum_{\alpha} \coh^{1,1}_{n,m,\alpha} \, {\bf 1}_{(m,n,\alpha)\in \mathcal{R}(k)}
\, e^{i \alpha_0 \cdot \theta_0},
{\bf N}^{1,0}_0+{\bf N}^{1,0}_{\rm space}
\right)\Big|_{T=t=0, \theta_1=0},
\end{align*}
an explicit, linear function of  $\underline{\bf U}^0$ and $\underline{\bf U}^1$.
Fourier transforming $\underline{\bf U}^0$ and $\underline{\bf U}^1$
in the variable $\theta_0$, and Fourier transforming
$$\left(
{\bf u}^{1,0}_0+\Pi {\bf u}^{1,0}_{\rm osc},
\sum_\alpha \coh^{1,1}_{n,m,\alpha} \, {\bf 1}_{(m,n,\alpha)\in \mathcal{R}(k)}
\, e^{i (\alpha_0 \cdot \theta_0+\alpha_1\cdot\theta_1))},
{\bf N}^{1,0}_0+{\bf N}^{1,0}_{\rm space}
\right)$$
in the variable $(\theta_0,\theta_1)$ then allows to deduce, as we did in the proof of
Theorem~\ref{ThExistPro}, the value of
$\left({\bf u}^{1,0}_0,
\Pi {\bf u}^{1,0}_\alpha,
\coh^{1,1}_{m,n,\alpha'},
{\bf N}^{1,0}_{\alpha''}
\right)\Big|_{T=t=0},
$
whenever
$\alpha\in\Z^{2d}\setminus\{0\}$,
$(m,n,\alpha')\in \mathcal{R}(k)$, 
and $\alpha''\in\cz\cup\{0\}$. The requirement \eqref{bibibi} then gives the value of
$\left({\bf u}^{1,0}_0,
\Pi {\bf u}^{1,0}_\alpha,
\coh^{1,1}_{m,n,\alpha'},
{\bf N}^{1,0}_{\alpha''}
\right)
$
on the whole set $\{T=0\}\cup\{t=0\}$.
Equation \eqref{intermtranspPop10}
in turn provides the value of
$\left({\bf u}^{1,0}_0,
\Pi {\bf u}^{1,0}_\alpha,
\coh^{1,1}_{m,n,\alpha'},
{\bf N}^{1,0}_{\alpha''}
\right)
$
for any value of $T$ and $t$.
The remaining part of ${\bf U}^1$, namely
the value of
$
\left((1-\Pi) {\bf u}^{1,0}_{\rm osc},
{\bf u}^{1,1}, 
{\bf C}^{1,0},
{\bf C}^{1,1}_{n,m,\alpha} \, {\bf 1}_{(m,n,\alpha)\notin \mathcal{R}(k)},
{\bf N}^{1,1}\right)
$
for all values of $T$ and $t$, is next given by $\phi\left({\bf U}^0\right)$.
This terminates the computation of ${\bf U}^1$ for all values of $T$ and $t$.
The relation 
${\bf U}^2=\psi\left({\bf U}^0,{\bf U}^1\right)$, see \eqref{boubouboubou},
prescribes 
${\bf U}^2$ for al values of $T$ and $t$.

Let us come to  regularity issues.
Thanks to the constraint \eqref{bouboubou},
the $\mathcal{C}^\ell(H^{s-\tilde{s}(\ell,a)})$ smoothness of 
$(1-\Pi) {\bf u}^{1,0}_{\rm osc}$,
${\bf u}^{1,1}$, 
${\bf C}^{1,0}$,
${\bf C}^{1,1}_{m,n,\alpha} \, {\bf 1}_{(m,n,\alpha)\notin \mathcal{R}(k)}$,
${\bf N}^{1,1}$ is simply the consequence of the $\mathcal{C}^\ell(H^{s-\tilde{s}(\ell,a)})$
smoothness of ${\bf U}^0$.
For 
$\fiel^{1,0}_0$, $\Pi\fuzo$, $\pop^{1,0}_\alpha \, {\bf 1}_{\alpha\in\cz\cup\{0\}}$ and 
$\coh^{1,1}_{m,n,\alpha} \, {\bf 1}_{(m,n,\alpha)\in\res}$,
the $\mathcal{C}^0(H^{s-\tilde{s}(0,a)})$ regularity  stems from the 
classical smoothness of solutions to linear hyperbolic systems with
source terms lying in $\mathcal{C}^0(H^{s-\tilde{s}(0,a)})$.
Note that in the last equation in \eqref{intermtranspPop10}, the regularity of the source term
$\dt \czu$ 
is a simple consequence of the regularity of solutions to ODE's depending on a parameter (here $t$).
Note also that the overall  loss
$\tilde{s}(0,a)$ comes from both the involved derivatives of ${\bf U}^0$
appearing in the source terms, and from
small divisors $(\omega(m,n)-\vf\cdot\alpha_1)^{-1}$ and 
$\muai$ acting on ${\bf U}^0$, that appear in the source term as well.
Differentiating \eqref{intermtranspPop10'} and \eqref{bouboubou}
with respect to $T$ and $t$, and applying the same argument,
eventually provides the 
$\mathcal{C}^\ell(H^{s-\tilde{s}(\ell,a)})$  smoothness of ${\bf U}^1$.
Relation ${\bf U}^2=\psi\left({\bf U}^0,{\bf U}^1\right)$ provides the
$\mathcal{C}^\ell(H^{s-\tilde{s}(\ell,a)})$  smoothness of ${\bf U}^2$.
\fin

\begin{rem}
One can prove that for any $\ell$
we have $\sol^1\in\mathcal{C}^\ell(H^{s-\max(3a+1,a+2)-\ell(a+2)})$, and
$\sol^2\in\mathcal{C}^\ell(H^{s-\max(4a+2,3a+3)-\ell(a+2)})$.
\end{rem}

\section{Convergence}
\label{SecCV}

\subsection{The residual} \label{Secres}

Writing down the full profile representation of the residual
$$r^\ep\var = \Lep\sola - F^\ep(\sola),$$
we have, using Proposition~\ref{PropDA},
\begin{align*}
& r^\ep\var = \mathcal{R}^\ep(t,x,y,z,T,\sigma,\theta)\Bigg|_{
T={t}/{\rep},
\sigma={\gamma t}/{\ep},(\theta_0,\theta_1)=\left(\vo x , -\vf t\right)/\ep},
\\
&
\text{with } \, 
\mathcal{R}^\ep = \sum_{j=-2}^3 \sum_{\kappa\geq0}\rep^j 
r^{j,\kappa}\varTt\exp(-\kappa\sigma).
\end{align*}
We now find how $r^1$, $r^2$ and $r^3$
depend on the profiles $\sol^1$,  $\sol^2$ and $\sol^3$. 
Symbolically, introducing first order differential operators $L_1$ and $L_2$, a matrix
$L_0$, and some bounded, symmetric, bilinear operators $B_E$, $B_{({\bf C}, {\bf N})}$
(with values on the $E$- and $({\bf C}, {\bf N})$-components, respectively), 
we have
\begin{equation*}
\begin{split}
r^1 & = L_2(\dT,\dy,\dz)\sol^2 + L_1(\dt,\dx)\sol^1 
+ 2 B_E(\sol^0,\sol^1) + 2 B_{({\bf C}, {\bf N})}(\sol^0,\sol^2) 
+ B_{({\bf C}, {\bf N})}(\sol^1,\sol^1) + L_0 \sol^1, \\
r^2 & = L_1(\dt,\dx)\sol^2 + 2 B_E(\sol^0,\sol^2) 
+ B_E(\sol^1,\sol^1) + 2 B_{({\bf C}, {\bf N})}(\sol^1,\sol^2) + L_0 \sol^2, \\
r^3 & = 2 B_E(\sol^1,\sol^2) + B_{({\bf C}, {\bf N})}(\sol^2,\sol^2).
\end{split}
\end{equation*}
Naturally, once $\mathcal{C}^\infty$ profiles have been built thanks to 
Theorem~\ref{ThExistPro} and Proposition~\ref{PropExistCorr}, 
regularity of these residual profiles is clear. 
The question we now need to face is to evaluate the size of these residuals
as the intermediate time $T$ grows unboundedly.

Concerning the dominant profile ${\bf U}^0$,
the components  $\fiel^0$ and $\pop^0$ are bounded uniformly in $T$,
while the component 
$\coh^0$ has exponential 
growth, {\em i.e.} $\coh^0$ and all its derivatives have size $K_1 \, \exp(K_2 \, T)$ as $T$ increases,
see Lemma~\ref{LemExp}.

Concerning the first corrector ${\bf U}^1$,
the components $\fiel^{1,0}_0$, $\Pi\fuzo$ and ${\bf N}^{1,0}_0+{\bf N}^{1,0}_{\rm space}$
are precisely constructed so as to be sublinear in $T$, see 
Section~\ref{SecInterm}. This is the key point.
On top of that,
the component $(1-\Pi) {\bf u}^{1,0}_{\rm osc}$
is bounded, thanks to
\eqref{polarFiel10},
the component ${\bf u}^{1,1}$ has exponential growth, thanks to
\eqref{polarFiel11},
the component ${\bf N}^{1,1}$ has exponential growth, thanks to
\eqref{polarPop1},
the component ${\bf C}^{1,0}$ is bounded, thanks to
\eqref{polarCoh10},
and the component ${\bf C}^{1,1}$has exponential growth, thanks to
\eqref{polarCoh11} and \eqref{intermtranspCoh11} in conjunction with the Gronwall Lemma.

Concerning the second corrector ${\bf U}^2$,
we know from \eqref{boubouboubou} that ${\bf U}^2$ is an explicit, linear function
of ${\bf U}^0$ and ${\bf U}^1$, functions that have at most exponential growth. Hence
${\bf U}^2$ has at most exponential
growth in $T$. The only difficulty may then come from
the component  ${\bf U}^{2,0}$, whose exponential growth will not  be eventually compensated
by a decaying term $e^{-\sigma}$ or so. In that direction, we observe
that ${\bf u}^{2,0}_0$, $\Pi {\bf u}^{2,0}_{\rm osc}$, and
${\bf N}^{2,0}_\alpha \, {\bf 1}_{\alpha\in\cz\cup\{0\}}$ are conventionally chosen to vanish.
On the other hand, the component
$(1-\Pi) {\bf u}^{2,0}_{\rm osc}$ is bounded thanks to
\eqref{polarFiel20} and to the boundedness of $({\bf u}^0,{\bf N}^0,{\bf u}^{1,0}_{\rm osc})$, 
the component ${\bf C}^{2,0}$
is bounded
thanks to \eqref{polarCoh21}, to the boundedness of ${\bf C}^{1,0}$, ${\bf u}^0$,  ${\bf N}^0$,
and to the sublinearity of  ${\bf N}^{1,0}$,
while the component
${\bf N}^{2,0}_\alpha \, {\bf 1}_{\alpha\notin\cz\cup\{0\}}$ 
is bounded
thanks to \eqref{polarPop20}, to the boundedness of ${\bf C}^{1,0}$, ${\bf u}^0$,  ${\bf N}^0$,
and to the sublinearity of  ${\bf N}^{1,0}$.

\medskip

As a conclusion, we have now established the
\begin{lemme} \label{LemEstimRes}
Given the $\mathcal{C}^\infty$ profiles provided by Theorem~\ref{ThExistPro} 
and Proposition~\ref{PropExistCorr} with the choice $s=+\infty$,
the following result holds.

For all $\mu\in\N^{5+2d}$,
there are constants $K_1,K_2>0$ such that, uniformly on 
$[0,t_\star]_t\times\R_T\times\R^3_{x,y,z}\times\T^{2d}_\theta$, 
for all $\kappa\in\N^\star$, $j=1,2,3$, we have
$$\left|\partial_{t,T,x,y,z,\theta}^\mu \, r^{j,\kappa}(T)\right|
\leq K_1e^{K_2T}, \,\mbox{ and }\, 
\frac{1}{T}\left|\partial_{t,T,x,y,z,\theta}^\mu \, r^{j,0}(T)\right| 
\tendlorsque{T}{+\infty} 0.$$
As a consequence, we get for the residual
$$\forall\mu\in\N^{3+d}, \quad \sup_{t\in[0,t_\star]} 
\Big\| \partial_{x,y,z,\theta_0}^\mu 
\mathcal{R}^\ep |_{T=t/\rep,\sigma=\gamma t/\ep,\theta_1=-\vf t/\ep} 
 \Big\|_{L^2_{x,y,z}} \tendlorsque{\ep}{0} 0,$$
\end{lemme}

\subsection{Stability} \label{SecStab}

In this section, we prove our main Theorem.
\begin{theo} \label{ThStab}
Let the profiles $\sol^0$, $\sol^1$, $\sol^2 \in 
\mathcal{C}^\infty([0,t_\star]\times[0,+\infty[^2,H^\infty
(\R^3\times\T^{2d}))$ be given by Theorem~\ref{ThExistPro} 
and Proposition~\ref{PropExistCorr}. They provide us with the 
approximate solution $\sola=\sola(t,x,y,z)$ given by equation \eqref{solapp}. 

Then, for any $s>(3+2d)/2$ and any familly  $(\dep)_{\ep>0} 
\subset H^s(\R^3\times\T^d)$ such that 
$\|\delta^\ep\|_{H^s} \tendlorsque{\ep}{0}0$, there is $\ep_0>0$ 
such that for $\ep\in]0,\ep_0]$, the Cauchy problem 
\begin{align}
\left\{
\begin{array}{l}
\vspace{0.3cm}
\displaystyle\Lep \sole = F^\ep(\sole) ,\\
\vspace{0.1cm}
 \sole_{|_{t=0}} = {\sola}_{|_{t=0}} + \dep(x,y,z,\vo x/\ep),
\end{array}
\right.
\end{align}
for the Maxwell-Bloch system \eqref{MBeps} has a unique (mild) solution 
$\sole$ which belongs to $\mathcal{C}([0,t_\star],H^s(\R^3))$. 
Besides, for all $\mu\in\N^3$ such that $s-|\mu|>(3+2d)/2$, we have
$$\|\partial^\mu_{x,y,z}\left(\sole-\sola\right)\|
_{L^\infty([0,t_\star]\times\R^3)} \tendlorsque{\ep}{0} 0.$$
\end{theo}

\noindent
{\bf Proof.}
Standard results for symmetric hyperbolic systems ensure that, for
$\ep>0$ fixed, a unique mild solution
$\sole\in\mathcal{C}([0,t_\ep],H^s(\R^3))$ exists for some $t_\ep>0$. 
The difficulty lies in bounding $t_\ep$ from below.

We use a singular system method (\cf \cite{JMR95}), and look for 
$\sole$ under the form of a profile, namely look for $\solc$ such that 
$$\sole\var=\solc(t,x,y,z,\vo x/\ep),$$
where $\solc=\solc(t,x,y,z,\theta_0)$ corresponds to the following 
initial data, which is non-singular in $\ep$ (this is a key point)
\begin{equation*}
\solc(0,x,y,z,\theta_0) = 
\left(
\sum_{j=0}^2 \rep^j 
\sol^j(t,x,y,z,T,\sigma,\theta_0,\theta_1) 
+ \dep(t,x,y,z,\theta_0) \right)_{|_{t=T=\sigma=0, \theta_1=0}}.
\end{equation*}
It is then sufficient, for $\sole$ to be a solution to \eqref{MBeps}, 
that $\solc$ satisfies
$$\Lsing\solc = F^\ep(\solc).$$
To go on with the analysis, we now set
\begin{align*}
& \solac(t,x,y,z,\theta_0) = \sum_{j=0}^2 \rep^j U^j
\varTdt_{|_{T=t/\rep,\sigma=\gamma t/\ep,\theta_1=\vf t/\ep}}.
\end{align*}
Note in passing that the at-most-exponential growth in $T$ of 
the various profiles $\sol^{j,\kappa}$ implies that the family 
$(\solac)_\ep$ is bounded in 
$\mathcal{C}([0,t_\star],H^s(\R^3\times\T^d))$. 
We evaluate the difference
\begin{align*}
& \Dep=\solc-\solac.
\end{align*}
In the next few lines, we may sometimes write
$\Dep=\left(\Delta_{\bf u}^\ep,\Delta_{\bf C}^\ep,\Delta_{\bf N}^\ep\right)$, refering to the ${\bf u}$, ${\bf C}$ and ${\bf N}$ components of $\Dep$, respectively. In any circumstance, 
we have $\Dep\in\mathcal{C}([0,t_\ep],H^s(\R^3\times\T^d))$, its initial value is
$\Dep_{|_{t=0}} = \dep$, and $\Dep$ satisfies
\begin{equation}
\label{systsing}
\Lsing \Dep = 
 F^\ep(\solac+\Dep)-F^\ep(\solac) -\mathcal{R}^\ep_{|_{T=t/\rep,\sigma=\gamma t/\ep,\theta_1=\vf t/\ep}}.
\end{equation}
Hence standard $H^s$ estimates provide
\begin{align*}
\frac{1}{2} \frac{\rm d}{\rm dt} \|\Dep\|_{H^s}^2 \leq 
\langle F^\ep(\solac+\Dep)-F^\ep(\solac) | \Dep\rangle_{H^s}
+
\left\|\mathcal{R}^\ep \right\|_{H^s} \, \left\|\Dep\right\|_{H^s}.
\end{align*}
Now, the evaluation of the scalar product
$\langle F^\ep(\solac+\Dep)-F^\ep(\solac) | \Dep\rangle_{H^s}$
involves various terms that may be ordered in powers of $1/\rep$
 as seen by inspection of the bilinear function $F^\ep$ in 
\eqref{MBeps}. The term carrying the weight $1/\ep$ is
\begin{align*}
&
\frac{1}{\ep}\langle -i \left[\Omega,\Delta_{\bf C}^\ep\right] -\gamma \, \Delta_{\bf C}^\ep \, | \, \Dep\rangle_{H^s}
=
-\frac{\gamma}{\ep}\langle   \Delta_{\bf C}^\ep \, |  \, \Dep\rangle_{H^s}
=
-\frac{\gamma}{\ep}\left\| \Delta_{\bf C}^\ep \right\|^2_{H^s},
\end{align*}
where the first equality comes from the fact that the operator $[\Omega , . ]$ is skew-symmetric.
To estimate the other terms, it is useful to 
keep in mind that $\solac$ is uniformly bounded in $H^s$ on the interval $[0,t_*]$, so there is a constant
$C$, independent of $f$ and $\ep$, such that  for any $t$ and $\ep$ we have $\left\|\solac(t)\right\|_{H^s}\leq C$.
Using this piece of information, the term carrying the weight $1$ in
$\langle F^\ep(\solac+\Dep)-F^\ep(\solac) | \Dep\rangle_{H^s}$
is clearly bounded by 
\begin{align*}
C \, \left\|\Dep\right\|_{H^s}^2,
\end{align*}
for some constant $C$ independent of $t$ and $\ep$, while the term carrying the weight $\rep$ in
is clearly bounded by 
\begin{align*}
\rep \,
\left\|\Dep\right\|_{H^s}^2,
\end{align*}
The more difficult term is the one carrying the weight $1/\rep$. Carefully treating apart
all occurences of the two terms $\Delta_{\bf C}^\ep$ and ${\bf C}^\ep_{\rm app}$,  all terms being majorized in the most simple fashion,
we recover that this contribution is upper-bounded by
\begin{equation*}
\frac{C}{\rep} \, \left(
\left\|\Delta_{\bf C}^\ep\right\|^2_{H^s}
+
\left\|\Delta_{\bf C}^\ep\right\|_{H^s}
\, \left\|\Delta^\ep\right\|_{H^s}
+
\left\|{\bf C}^\ep_{\rm app}\right\|_{H^s}
\, \left\|\Delta^\ep\right\|_{H^s}^2  
+
\left\|{\bf C}^\ep_{\rm app}\right\|_{H^s}
\, \left\|\Delta_{\bf C}^\ep\right\|_{H^s}
\, \left\|\Delta^\ep\right\|_{H^s}  
\right).
\end{equation*}
All in all we have eventually proved, gathering some terms for convenience,
\begin{align*}
&
\frac{\rm d}{\rm dt} \|\Dep\|_{H^s}^2
\leq 
-\frac{\gamma}{\ep}\left\| \Delta_{\bf C}^\ep \right\|^2_{H^s}
+
\frac{C}{\rep} \,  \left( \left\|\Delta_{\bf C}^\ep\right\|^2_{H^s}
+
\left\|\Delta_{\bf C}^\ep\right\|_{H^s} \, \left\|\Delta^\ep\right\|_{H^s}\right)
\\
&
\qquad \qquad \qquad+
\frac{C}{\rep} \, \left\|{\bf C}^\ep_{\rm app}\right\|_{H^s} \, 
\left(
\left\|\Delta_{\bf C}^\ep\right\|_{H^s}
\left\|\Delta^\ep\right\|_{H^s}
+
\left\|\Delta^\ep\right\|_{H^s}^2  
\right)
\\
&
\qquad \qquad \qquad
+
C \, \left\|\Dep\right\|_{H^s}^2
+
\rep \,
\left\|\Dep\right\|_{H^s}^2
+
\left\|\mathcal{R}^\ep \right\|_{H^s} \, \left\|\Dep\right\|_{H^s}.
\end{align*}
Hence, there is an $\ep_0$ such that for any $\ep\leq\ep_0$ we have, 
\begin{equation}
\label{grongron}
\frac{\rm d}{\rm dt} \|\Dep\|_{H^s}^2 \leq 
-\frac{1}{C \, \ep} \, \left\| \Delta_{\bf C}^\ep \right\|^2_{H^s}
+
\frac{C}{\rep} \, \left\| \Delta_{\bf C}^\ep \right\|_{H^s} \, \left\| \Delta^\ep \right\|_{H^s}
+
C \, \frac{\left\|{\bf C}^\ep_{\rm app}\right\|_{H^s} }{\rep} \, \left\| \Delta^\ep \right\|_{H^s}^2
+
\left\|\mathcal{R}^\ep \right\|_{H^s} 
\,
\left\|\Dep\right\|_{H^s}.
\end{equation}
Now, the two crucial ingredients are
\begin{align}
\label{grou}
&
\mathop{\rm sup}_{t\in[0,t_*]} \left\|\mathcal{R}^\ep \right\|_{H^s} 
\mathop{\to}_{\ep \to 0} 0,
\end{align}
thanks to Lemma \ref{LemEstimRes} , and
\begin{align}
\label{tougrou}
\int_0^{t_*} \frac{\|\coha(t)\|_{H^s}}{\rep} {\rm d}t \leq C \, \rep,
\end{align}
for some $C$ independent of $\ep$.
This crucial piece of information comes from the fact that
\begin{align*}
\coha(t)={\bf C}^{0,1}\big|_{T=t/\rep} \, e^{-\gamma t/\ep} + \rep \, {\bf C}^{1,1}\big|_{T=t/\rep} \, e^{-\gamma t/\ep} 
+
\mathcal{O}(\ep),
\end{align*}
which, in conjunction with the at-most-exponential growth of ${\bf C}^{0,1}$ and of ${\bf C}^{1,1}$ (see {\em e.g.}
Lemma \ref{LemExp} and
estimate \eqref{estimCoh11}), provides
\begin{align*}
\|\coha(t)\|_{H^s}
\leq
\exp\left(-\frac{t}{C \, \ep}\right)+ C \, \ep,
\end{align*}
for some $C$ independent of $t\in[0,t_*]$ and $\ep$.

At this stage, an easy argument using the Gronwall Lemma allows to deduce from
\eqref{grongron}, \eqref{grou}, and \eqref{tougrou}, that $\|\Dep\|_{H^s}$ is bounded independently of $\ep$ over the whole
interval $[0,t_*]$. Repeting the Gronwall argument next shows that  $\|\Dep\|_{H^s}$ actually satisfies
\begin{equation}
\label{grongrongron}
\|\Dep\|_{H^s} \leq 
\|\delta^\ep\|_{H^s}
\, \exp\left( 
C \displaystyle \int_0^t  \left(
\frac{\left\|{\bf C}^\ep_{\rm app}\right\|_{H^s}}{\rep} 
+
\left\|\mathcal{R}^\ep \right\|_{H^s} 
\right) \, dt'
\right)
\mathop{\longrightarrow}_{\ep\to 0} 0, \quad
\text{uniformly on } [0,t_*].
\end{equation}
We skip the whole Gronwall-like argument.
Estimate \eqref{grongrongron} now induces, by Sobolev's injection, $L^\infty$ convergence to zero  for 
profiles $\Dep=\solc-\solac$, and thus $L^\infty$ convergence to zero for the original functions $\sol-\sola$, since 
$\|\sol-\sola\|_{L^\infty_{t,x,y,z}}=
\|\solc-\solac\|_{L^\infty_{t,x,y,z,\theta_0}}$.
\fin

\section{The Transverse Magnetic case} \label{SecTM}

\subsection{The system}

In this section, we present the previous WKB in the particular
Transverse Magnetic case, when fields take the form
$$
{B}=
\begin{pmatrix}
B_x \\
B_y \\
0
\end{pmatrix}
={B}(t,x,y),
\quad
{E}=
\begin{pmatrix}
0 \\
0 \\
E
\end{pmatrix}
={E}(t,x,y),
$$
with the additional common assumption (\cf \cite{NM92}) that the 
polarization operator $\Gamma$ has entries parallel to ${E}$, namely
$$\forall m,n \in \{1,\dots,N\}, \quad{\Gamma}(m,n)=
\begin{pmatrix}
0 \\
0 \\
\Gamma(m,n)
\end{pmatrix}.
$$
Maxwell-Bloch system then reads 
\begin{align}
\nonumber
&
\dt B_x^\ep + \frac{1}{\rep}\dy E^\ep = 0,
\\
\nonumber
&
\dt B_y^\ep + \frac{1}{\rep}\dx E^\ep = 0,
\\
&
\nonumber
\dt E^\ep - \dx B_y^\ep  + \dy B_x^\ep = 
\frac{i}{\rep} \Tr \left(\Gamma\Og\cohe\right) 
- i E^\ep \Tr \left(\Gamma[\Gamma,\cohe+\pope]\right) 
- \rep \Tr \left(\Gamma W\op \pope\right),
\\
\label{TMMBeps}
&
\dt\cohe = -\frac{i}{\ep} \Og \cohe 
+ \frac{i}{\rep} E^\ep \crod{\Gamma}{\cohe+\pope} ,
\\
&
\nonumber
\dt\pope = \frac{i}{\rep}E^\ep\crd{\Gamma}{\cohe} + W\op\pope.
\end{align} 

\subsection{The Ansatz}

As stressed in Remark~\ref{RemAnsatz}, the introduction of an
intermediate time $T=t/\rep$ is not necessary here (see also 
Remark~\ref{noT} below). In order to simplify computations, 
we also restrict here our attention  to the case of prepared data, which corresponds to the case when
$$
{\bf C}^\ep\big|_{t=0}=0.
$$
This second simplification allows us not to
 use 
the variable $\sigma=\gamma t/\ep$ in the sequel. In a nutshell, we here consider the simplified Ansatz
\begin{equation}
\label{TMsolapp}
\sola(t,x,y) = \sum_{j=0}^2 \rep^j\sol^j(t,x,y,\theta)
_{|_{\theta=(\vo x/\ep,-\vf t/\ep}}, \quad 
\sol^j(t,x,y,\theta)
=
\sum_{\alpha\in\Z^{2d}} \sum_{\kappa\in\N}
\sol^j_\alpha(t,x,y) \,
e^{i\alpha\cdot\theta} .
\end{equation}

The characteristic sets $\cz$, $\cp$, $\cm$ and the resonant set 
$\res$ are the same as before.
\subsection{WKB expansions}

In this setting, vanishing of the terms $r^{-2}$, $r^{-1}$ and 
$r^0$ from Proposition~\ref{PropDA} reduces to
\vspace{-0.1cm}
\begin{align} \label{TMres-2}
\left\{
\begin{array}{l}
\vspace{0.15cm}
-\vf\cdot\dthu B_x^0 = 0 , \\
\vspace{0.15cm}
 -\vf\cdot\dthu B_y^0 - \vo\cdot\dthz E^0 = 0 , \\
\vspace{0.15cm}
 -\vf\cdot\dthu E^0 - \vo\cdot\dthz B_y^0 = 0 , \\
\vspace{0.15cm}
 (i\,\Og-\vf\cdot\dthu)\coh^0 = 0 , \\
 -\vf\cdot\dthu\pop^0 = 0 , \\
\end{array}
\right.
\end{align}
\vspace{-0.4cm}
\begin{align}\label{TMres-1}
\left\{
\begin{array}{l}
\vspace{0.15cm}
 -\vf\cdot\dthu B_x^1 + \dy E^0 = 0 , \\
\vspace{0.15cm}
 -\vf\cdot\dthu B_y^1 - \vo\cdot\dthz E^1 = 0 , \\
\vspace{0.15cm}
 -\vf\cdot\dthu E^1 - \vo\cdot\dthz B_y^1 +\dy B_x^0 = 0 , \\
\vspace{0.15cm}
 (i\,\Og-\vf\cdot\dthu)\coh^1  
= i\, E^0 \crod{\Gamma}{\coh^0+\pop^0} , \\ 
 -\vf\cdot\dthu\pop^1 
= i\, E^0 \crd{\Gamma}{\coh^0} ,
\end{array}
\right.
\end{align}
\begin{align} \label{TMres0}
\left\{
\begin{array}{l} 
\vspace{0.15cm}
-\vf\cdot\dthu B_x^2 + \dy E^1 + \dt B_x^0 = 0 , \\
\vspace{0.15cm}
 -\vf\cdot\dthu B_y^2 - \vo\cdot\dthz E^2 + \dt B_y^0 
- \dx E^0 = 0 , \\
\vspace{0.15cm}
 -\vf\cdot\dthu E^2 - \vo\cdot\dthz B_y^2 +\dy B_x^1 
+ \dt E^0 - \dx B_y^0 
= i\, \Tr(\Gamma\Og\coh^1) 
-i\, E^0 \Tr(\Gamma[\Gamma,\coh^0+\pop^0]) , \\
\vspace{0.15cm}
 (i\,\Og-\vf\cdot\dthu)\coh^2  + \dt \coh^{0,\kappa} 
= i\, E^0 \crod{\Gamma}{\coh^1+\pop^1} 
+ i\, E^1 \crod{\Gamma}{\coh^0+\pop^0} ,  \\
  -\vf\cdot\dthu\pop^2 + \dt \pop^0 
= i\, E^0 \crd{\Gamma}{\coh^1} 
+ i\, E^1 \crd{\Gamma}{\coh^0} + W\op \pop^0 .
\end{array}
\right.
\end{align}

According to \eqref{TMres-2}, the polarization conditions from 
equations~\eqref{polarFiel}--\eqref{polarPop} become, with the 
notations from Definition~\ref{DefOsc}
\begin{align}
&
\label{TMpolarFielLong}
B_{x,{\rm time}}^0 = 0 ,
\\
&
\label{TMpolarFielTrans1}
\forall \alpha \notin \cp\cup\cm\cup\{0\}, 
\quad B_{y,\alpha}^0 = E^0_\alpha = 0 ,
\\
&
\label{TMpolarFielTrans2}
\forall \alpha \in \cp\cup\cm, 
\quad B_{y,\alpha}^0 = \mp E^0_\alpha  ,
\\
&
\label{TMpolarCoh}
\coh^0 = 0 ,
\\
&
\label{TMpolarPop}
\pop^0_{\rm time} = 0 .
\end{align}
Here, the projector $\Pi$ from Definition~\ref{DefProj} is 
given explicitly, and the ``prepared data'' condition is coherent with the 
vanishing of $\coh^0$.

From \eqref{TMres-1}, we get for the average $\fiel_0^0$ a 
version of the evolution equation at intermediate scale 
\eqref{intermtranspFielMoy0} under the form
$$\dy E_0^0 = \dy B_{x,0}^0 = 0 ,$$
which leads to  polarization consitions, 
instead of evolution equations, namely
\begin{equation} \label{TMpolarFielMoy}
E_0^0 = B_{x,0}^0 = 0 .
\end{equation}
The oscillating part of \eqref{TMres-1} produces the following
transcription of the polarization condition \eqref{polarFiel10} 
for the first corrector $\fiel^1$
\begin{align} \label{TMpolarFiel1}
\left\{
\begin{array}{l}
\vspace{0.2cm}
 \forall \alpha \notin \ch\cup\{0\}, 
\quad \fiel^1_\alpha = 0 , \\
\vspace{0.2cm}
 \forall \alpha \in \cz, \quad E^1_\alpha=0, \quad 
B_{y,\alpha}^1 = \frac{1}{i\vo\cdot\alpha_0} \dy B_{x,\alpha}^0 , \\
 \forall \alpha \in \ch_\pm, \quad 
B_{x,\alpha}^1 = \frac{1}{i\vf\cdot\alpha_1} \dy E_\alpha^0 ,
\quad B_{y,\alpha}^1 = \mp E_{\alpha}^1 .\\
\end{array}
\right.
\end{align}
For coherences and populations, equation~\eqref{polarCoh10} is
unchanged,
\begin{equation} \label{TMpolarCoh1}
\coh^1 = i\,(i\,\Og-\vf\cdot\dthu)^{-1} E^0 \crd{\Gamma}{\pop^0},
\end{equation}
and equation~\eqref{polarPop1} becomes
\begin{equation} \label{TMpolarPop1}
\pop^1_{\rm time} =0 .
\end{equation}
Next, the average of fields equations in \eqref{TMres0} is 
equivalent to
\begin{equation} \label{TMpolarEMoy1}
\dy E^1_0 = 0, \qquad \ie \qquad E^1_0 = 0 
\end{equation}
(a polarization condition, again, instead of evolution as in 
equation~\eqref{intermtranspFielMoy10}),
\begin{equation} \label{TMslowtranspBy0}
\dt B_{y,0}^0 = 0 ,
\end{equation}
(playing the role of equation~\eqref{slowtranspFielMoy0}), 
and 
\begin{equation*}
\dy B_{x,0}^1 - \dx  B_{y,0}^0 
= - i (E^0 \Tr(\Gamma[\Gamma,\pop^0]))_0 ,
\end{equation*}
which reduces to
\begin{equation} \label{TMpolarBMoy1}
\dy B_{x,0}^1 - \dx  B_{y,0}^0 = 0 ,
\end{equation}
thanks to polarizations conditions and spectral properties \eqref{TMpolarFielTrans1}, 
\eqref{TMpolarPop} and \eqref{TMpolarFielMoy}.

\begin{rem} \label{noT}
Since $B_{y,0}^0$ does not depend on time $t$ (according to 
equation~\eqref{TMslowtranspBy0}), we may impose (as in 
equation~\eqref{TMpolarBMoy1}) that $\dx  B_{y,0}^0$ be the 
$y$-derivative of an $H^s$ function for all times, simply by 
requiring this condition be satisfied at $t=0$. But in the general three-dimensional 
framework, $\dt B_{y,0}^0 = 0$ is not given \emph{a priori}, 
and we need the addition of the intermediate variable $T$ to 
perform the analysis of Section~\ref{SecInterm}. Omitting this 
intermediate time leads, in the three-dimensional case, to the 
overdetermined (and ill-posed) system 
\eqref{intermtranspFielMoy0}, \eqref{intermtranspFielMoy10}
$$\mdyz\fiel^0_0=0 , \quad\mdyz\fiel^1_0=-\mun\fiel^0_0.$$
\end{rem}

Oscillations in \eqref{TMres0} are analyzed as follows. The 
polarization~\eqref{polarFiel20} for $\fdo$ splits into
\begin{equation} \label{TMpolarBx2}
B_{x,{\rm time}}^2 = (\vf\cdot\dthu)^{-1} \dy E^1_{\rm time} ,
\end{equation}
and
\begin{align} \label{TMpolarFiel2}
\left\{
\begin{array}{l}
\vspace{0.15cm}
 \forall \alpha \in \ch_\pm, 
\quad 2\, i \, (B_{y,\alpha}^2 \pm E^2_\alpha) =  
\displaystyle \frac{1}{i\vf\cdot\dthu} \dy^2 E^0_\alpha 
- i \Tr(\Gamma\Og\coh^1_\alpha) 
+ i (E^0\Tr(\Gamma\Og\pop^0)_\alpha , \\
\vspace{0.15cm}
 \forall \alpha \notin \cp\cup\cm\cup\{0\} , \quad
 \begin{pmatrix}
B_{y,\alpha}^2 \\
E^2_\alpha
\end{pmatrix}
= \displaystyle\frac{i}{(\vo\cdot\alpha_0)^2-(\vf\cdot\alpha_1)^2}
\times \\
\qquad\qquad\qquad\qquad\qquad\qquad \times 
\begin{pmatrix}
\vf\cdot\alpha_1 & -\vo\cdot\alpha_0 \\
-\vo\cdot\alpha_0 & \vf\cdot\alpha_1
\end{pmatrix} 
\begin{pmatrix}
0 \\
\dy B_{x,\alpha}^1 + i \Tr(\Gamma\Og(E^0\pop^0-\coh^1)_\alpha) 
\end{pmatrix},
\end{array}
\right.
\end{align}
where $B_{x,\alpha}^1$ vanishes for $\alpha\notin\ch$ (because 
of \eqref{TMpolarFiel1}), and may be chosen arbitrarily when 
$\alpha\in\cz$. The coherence $\coh^1$ is given by equation
~\eqref{TMpolarCoh1}. Then, the evolution equation
~\eqref{slowtranspFielOsc0} with respect to the slow time $t$
corresponds to
\begin{align} \label{TMslowtranspBx0}
\dt B_{x,{\rm space}}^0 = 0, \quad
\text{ together with } \quad 
\dt B_{x,0}^0 = 0,
\end{align}
and, using \eqref{TMpolarCoh1}, to
\begin{equation*}
\forall\alpha\in\ch_\pm, \quad
2(\dt\pm\dx)E^0_\alpha 
+ \frac{1}{i\vf\cdot\alpha_1}\dy^2 E^0_\alpha 
= i\,\Tr(\Gamma\Og\coh^1_\alpha)
-i\,(E^0\Tr(\Gamma[\Gamma,\pop^0]))_\alpha ,
\end{equation*}
or (with $\coh^1_\alpha$ from \eqref{TMpolarCoh1}), in other words
\begin{equation} \label{TMslowtranspE0}
2(\dt+\vDth\dx)E^0 
+ \dthu^{-1}\dy^2 E^0 \\
= i\,\Tr(\Gamma(i\,\Og(i\,\Og-\vf\cdot\dthu)^{-1}-1)
(E^0[\Gamma,\pop^0])),
\end{equation}
(in place of equation~\eqref{intermtranspFielOsc10}), with 
$\vDth$ given in Lemma~\ref{LemVg}.

For coherences, we get
\begin{equation} \label{TMpolarCoh2}
\coh^2=i\,(i\,\Og-\vf\cdot\dthu)^{-1}
(E^1\crod{\Gamma}{\pop^0}+E^0\crod{\Gamma}{\coh^1+\pop^1}),
\end{equation}
and for populations, we have 
\begin{equation} \label{TMpolarPop2}
\pop^2_{\rm time} = -i\,(\vf\cdot\dthu)^{-1}
(E^0\crod{\Gamma}{\coh^1})_{\rm time},
\end{equation}
and 
$$\dt\pop^0=W_\op\pop^0+i\,(E^0\crd{\Gamma}{\coh^1})_{\rm space}+i\,(E^0\crd{\Gamma}{\coh^1})_{0},$$
{\rm i.e.}
\begin{align} \label{TMslowtranspPop0}
&
\dt\pop^0 = W_\op\pop^0 
- (E^0\crd{\Gamma}{(i\,\Og-\vf\cdot\dthu)^{-1}
(E^0\crod{\Gamma}{\pop^0})})_{\rm space}\\
&
\nonumber
\qquad\qquad\qquad
- (E^0\crd{\Gamma}{(i\,\Og-\vf\cdot\dthu)^{-1}
(E^0\crod{\Gamma}{\pop^0})})_{0}.
\end{align}

\subsection{Conclusion in the TM case}

The above  computations provide us with a set of profile 
equations leading to a local in time, smooth approximate 
solution $\sola$ as in \eqref{TMsolapp}: 

$\bullet$ The leading profile $(\fiel^0,\coh^0,\pop^0)$ is 
given by polarizations \eqref{TMpolarFielLong}, 
\eqref{TMpolarFielTrans1}, \eqref{TMpolarFielTrans2}, 
\eqref{TMpolarCoh} ($\coh^0=0$), \eqref{TMpolarPop}, as well 
as \eqref{TMpolarFielMoy} (for average of fields, instead of an
evolution at the intermediate scale). They also satisfy  an 
evolution with respect to time $t$, given by the trivial equations 
\eqref{TMslowtranspBy0} and \eqref{TMslowtranspBx0}, namely
$$\dt B_{y,0}^0 = \dt B_{x,{\rm sp}}^0 = 0,$$
by a nonlinear Schr\"odinger equation~\eqref{TMslowtranspE0}, 
and by the Boltzmann equation~\eqref{TMslowtranspPop0}.
As quoted in Remark~\ref{noT}, this \emph{a priori} 
overdetermined set of equations is in fact well-posed because 
of its scalar structure. We stress the fact that no 
rectification occurs for fields at leading order: the only 
non-vanishing average $B^0_{y,0}$ is constant in time. \\
$\bullet$ Correctors are partially determined by polarizations, namely 
\eqref{TMpolarFiel1}, \eqref{TMpolarEMoy1}, 
\eqref{TMpolarBMoy1} for $\fiel^1$, 
\eqref{TMpolarCoh1} for $\coh^1$, 
\eqref{TMpolarPop1} for $\pop^1$, 
\eqref{TMpolarBx2}, \eqref{TMpolarFiel2} for $\fiel^2$, 
\eqref{TMpolarCoh2} for $\coh^2$, and 
\eqref{TMpolarPop2} for $\pop^2$. Parts of the correctors that are
not submitted to these constraints may be chosen equal to zero. 

The Ansatz
$$\sola(t,x,y) 
= \solac(t,x,y,\theta)_{|_{\theta=(\vo/\ep,-\vf t/\ep)}}$$  
is consistant with system~\eqref{TMMBeps}. We easily prove the
\begin{lemme} \label{LemEstimResTM}
Given the $\mathcal{C}^\infty$ profiles above, define the 
residual $\mathcal{R}^\ep$ as in Proposition~\ref{PropDA}. 
Then, for all $\mu\in\N^{2+d}$, there is $K>0$ such that:
$$\sup_{t\in[0,t_\star]} 
\Big\| \partial_{x,y,\theta_0}^\mu \mathcal{R}^\ep 
|_{\theta_1=-\vf t/\ep} \Big\|_{L^2_{x,y}} \leq K \rep .$$
\end{lemme}

Note that for prepared data, no initial layer is created at 
leading order, so that we get a $\mathcal{O}(\rep)$ estimate 
instead of the $o(1)$ in Lemma~\ref{LemEstimRes}.

By the same technique as in Section~\ref{SecStab}, we get 
finally 
\begin{theo} \label{ThStabTM}
Given the smooth profiles above on $[0,t_\star]$, for 
$s>(2+d)/2$ and any familly $(\dep)_{\ep>0} \subset 
H^s(\R^2\times\T^d)$ such that $\|\delta^\ep\|_{H^s} = 
\mathcal{O}(\rep)$, there is $\ep_0>0$ such that for 
$\ep\in]0,\ep_0]$, the Cauchy problem for Transverse Magnetic 
Maxwell-Bloch system \eqref{TMMBeps}, with initial data 
${\sola}_{|_{t=0}}+\delta^\ep(x,y,\vo x/\ep)$, has a unique 
solution $\sole\in\mathcal{C}([0,t_\star],H^s(\R^2))$, 
and for all $\mu\in\N^2$ such that $s-|\mu|>(2+d)/2$, 
there is $K>0$ such that 
$$\|\partial^\mu_{x,y}\left(\sole-\sola\right)\|
_{L^\infty([0,t_\star]\times\R^2)} \leq K \rep.$$
\end{theo}

\bibliographystyle{plain}

\end{document}